\documentclass[11pt]{article}

\usepackage{amsmath,amssymb}

\input xy \xyoption{all}
\CompileMatrices

\ifx\ThmNum\undefined 
   \newtheorem{theorem}{Theorem}[section]
   
\else 
   \newtheorem{theorem}{Theorem} 
   
\fi

\newtheorem{proposition}[theorem]{Proposition}
\newtheorem{lemma}[theorem]{Lemma}
\newtheorem{definition}[theorem]{Definition}
\newtheorem{conjecture}[theorem]{Conjecture}

\newtheorem{corollary}[theorem]{Corollary}

\def\proof{\par\medskip\noindent {\bf Proof: }}

\def\proofof #1 {\par\medskip\noindent {\sc Proof of #1. }}
\def\Box{\framebox[10pt]{\rule{0pt}{3pt}}}

\def\qed{\hfill $\Box$ \medskip \par}

\def\remark{\par\medskip \noindent {\bf Remark: }}

\def\map{{\bf map}}
\def\Hom{{\bf Hom}}

\def\A{\mathcal{A}}
\def\C{\mathcal{C}}

\def\E{\mathcal{E}}
\def\D{\mathcal{D}}
\def\G{\mathcal{G}}
\def\H{\mathcal{H}}
\def\I{\mathcal{I}}
\def\U{\mathcal{U}}
\def\M{\mathcal{M}}
\def\OO{\mathcal{O}}
\def\W{\mathcal{W}}

\newcommand{\Einf}{{E_{\infty}}}
\newcommand{\catM}{{\mathcal{M}}}
\def\F{\mathcal{F}}
\def\S{\Delta^{op}}
\def\SC{\Delta^{op}PrShv(\C)}
\def\SM{\Delta^{op}PrShv(\M)}
\def\SCp{\Delta^{op}PrShv(\C) _{\bullet}}
\def\Sp{Sp^{\Sigma}}
\def\SpC{Sp^{\Sigma,T}(\C)}
\def\SpM{Sp^{\Sigma,T}(\M)}
\def\Gm{{\bf G}_m}
\def\Z{\mathbf Z}
\def\Ab{\bf Ab}
\def\smash{\wedge}
\title{\vskip -25mm Preorientations of the derived motivic
multiplicative group}
\author{Jens Hornbostel}
\date{\today}

\begin{document}

\maketitle

\begin{abstract}
We establish several new model structures
and Quillen adjunctions both in the classical 
and in the motivic case for algebras over
operads and for modules over strictly
commutative ring spectra. As an application,
we provide a proof in the language of model categories
and symmetric spectra of Lurie's theorem that 
topological complex $K$-theory represents 
orientations of the derived multiplicative group.
Then we generalize this result to the motivic situation.
The appendix contains an erratum concerning this 
last point.
(AMS subject classification: Primary 55U35,
Secondary 18D10, 19D06, 55P48.) 

\end{abstract}

\section{Introduction}

{\bf Added in 2017}: {\em This version contains an appendix
with a short erratum written in 2017. As the preprint was published in 2013, 
we did not make any corrections in the preprint itself.}

\medskip

In this article, we establish several new model structures
and Quillen adjunctions in the motivc setting and study their
basic properties. In particular, we establish 
stable positive model structures for algebras over operads 
in motivic symmetric spectra. Moreover, we show that
both the flat and the projective stable positive
model structures on motivic symmetric spectra
satisfy the Goerss-Hopkins axioms, and that the flat variant
lifts to a model structure on strictly commutative
ring spectra which satisfies the To\"en-Vezzosi axioms.
In short, we make a first step towards a motivic version
of derived algebraic geometry.

\medskip

After the invention of several strict monoidal model
categories underlying the stable homotopy 
category (see \cite{EKMM}, \cite{HSS}), the subsequent
study of commutative algebra of strictly commutative
ring spectra has drawn
a lot of attention in recent years.
Glueing these ``derived'' commutative ring objects
together leads to one of the possible 
frameworks for derived algebraic geometry, with
classical algebraic geometry embedded via the 
Eilenberg-Mac Lane functor.
Derived algebraic geometry through commutative ring
spectra has become even more popular
when Jacob Lurie gave a conceptual definition
of $tmf$ (that is topological modular forms)
as the solution of a moduli problem in  
derived algebraic geometry, as opposed to the
handicraft construction of Goerss-Hopkins -Miller
\cite{Be}. More precisely, Lurie constructs
$tmf$ are the global sections of a sheaf of
$E_{\infty}$-ring spectra classifying
oriented derived elliptic curves.
He has sketched the proof of this theorem
in \cite{Lu}, and has lectured about various parts
of it at various places. His point of view
is that the best language to state and prove
the theorem is the one of infinity categories
rather than the one of model categories,
and we have no reason to doubt he is right.  
Here infinity categories really mean
quasi-categories, also known as weak
Kan complexes, as first invented by 
Boardman and Vogt \cite{BV} and recently studied in 
great detail by Lurie, Joyal and others.
The interested reader should consult
Lurie's homepage and \cite{Lu2}, as well as
\cite{Ber} for a comparison with other approaches
to infinity categories.
We expect that Lurie will publish a detailed proof 
of his theorem in this language in the near future,
the book \cite{Lu2} and the preprints of Lurie containing
already most of the necessary language and machinery.

\medskip

The above description of $tmf$ 
(corresponding to height 2 and the second 
chromatic layer) has an analog in
height 1 which is much easier to state
and to prove, and is also due
to Lurie \cite[section 3]{Lu}. Namely, real 
topological $K$-theory $KO$ classifies oriented derived
multiplicative groups. The key step
for proving this is to show that the suspension spectrum
of ${\bf CP}^{\infty}$ classifies
preorientations of the derived multiplicative group.
Here the derived multiplicative group is
by definition ${\bf G}_m:=Spec(\Sigma^{\infty}{\bf Z}_+)$,
the name being justified by
classical algebraic geometry over a base field $k$,
where the multiplicative group is $Spec(k[{\bf Z}])$.
As usual, the object $Rmap_{AbMon(\Sp)}(\Sigma^{\infty} \Z_+,
-))$ it represents
via the derived version of the 
Yoneda embedding will still be called
the multiplicative group.
(In the present preprint, all arguments
take place in the {\it affine} derived setting, 
so there is no need to write $Spec$ and to reverse 
the order of the arrows everywhere.)
We will provide a proof of this result
in the language of model categories and symmetric
spectra. For this we tacitly assume ${\bf CP}^{\infty}$
(compare \cite{Lu}) that
the topological monoid ${\bf CP}^{\infty}$ has been 
replaced by a homotopy equivalent model that
is actually a topological resp. simplicial group.
Then the result reads as follows
in general, the special case $N={\bf CP}^{\infty}$
being the one discussed above:
\begin{theorem}\label{main1}
(Lurie) For any abelian monoid $A$ in symmetric spectra
$\Sp$ (based on simplicial sets) and any simplicial abelian group 
$N$,
we have a natural isomorphism 
of abelian groups

$$ Hom_{Ho(AbMon(\Sp))}(\Sigma^{\infty}N_+,A)$$
$$ \simeq 
Hom_{Ho(AbMon(\S Sets))}(N,Rmap_{AbMon(\Sp)}(\Sigma^{\infty} \Z_+,
A))$$
$$ =
Hom_{Ho(AbMon(\S Sets))}(N,\Gm (A)).$$
\end{theorem}

Here $Ho(-)$ denotes the homotopy category,
$Rmap$ means the derived mapping space 
and the weak equivalences between abelian monoids
are always the underlying
ones, forgetting the abelian monoid structure.
The model structures involved in this statement are 
discussed in detail in section 3. Beware that in general
the category of abelian monoids in a homotopy category
of a monoidal model category is different 
from the homotopy category of abelian monoids in the 
monoidal model category, the monoidal model category
$(\S Sets, \times)$
and the abelian monoid $QS^0$ in $Ho(\S Sets)$
being the most prominent example.

\medskip 

We will prove Theorem \ref{main1} in section 4,
and indeed a much more general version
(see below). Using a theorem of Snaith \cite{Sn}, 
Lurie's definition of an orientation and his above
theorem then imply that
complex $K$-theory represents orientations
of the standard derived multiplicative group,
and then further that real $K$-theory represents
oriented derived multiplicative groups
in general.
We refer the reader to section 5 for further
details.

\medskip

The above theorem can be generalized 
to motivic symmetric spectra
$\SpM$ on the smooth Nisnevich site $\M=(Sm/S)_{Nis}$ with $S$
an arbitrary noetherian base scheme as follows, 
everything equipped with appropriate
motivic (that is ${\bf A}^1$-local)
model structures 
as discussed in section 3,
and the motivic derived multiplicative
group $\Gm^{mot}$ defined using
the suspension spectrum
with respect to a given motivic circle
$T$, that is represented by $\Sigma_T^{\infty}{\bf Z}_+$.
(As the notation suggests, we are mainly
interested in the circle provided by 
${\bf A}^1/({\bf A}-0)$ resp. the weakly
equivalent $T={\bf P}^1$ pointed at 
infinity.)

\begin{theorem}\label{main2}
Let $\M=(Sm/S)_{Nis}$ and $T=S^1$ or $T={\bf P}^1$. 
Then for any abelian monoid $A$ in motivic
symmetric $T$-spectra $\SpM$ and any abelian group 
$N$ in the category $\SM$ of 
simplicial presheaves on $\M$, we have a natural isomorphism 
of abelian groups

$$Hom_{Ho(AbMon(\SpM))}(\Sigma^{\infty}_{T}N_+,A)$$ 
$$\simeq Hom_{Ho(AbMon(\SM))}(N,Rmap_{AbMon(\SpM)}(\Sigma^{\infty}_{T}\Z_+,
A))$$
$$ = Hom_{Ho(AbMon(\SM))}(N,\Gm (A)).$$
\end{theorem}

This is a rather straightforward application of the main
technical results of this article.
Appliying it to $T={\bf P}^1$ pointed at ${\infty}$
and to $N={\bf P}^{\infty}$ which is 
not a variety but still a simplicial presheaf,
and using the recently established motivic version
of Snaith's theorem \cite{GS}, \cite{SO},
it will imply that {\it algebraic $K$-theory 
represents motivic orientations of the 
derived motivic multiplicative group},
provided one works with the correct motivic
generalizations of the concept of derived algebraic 
groups and of orientations.
Again, we refer to section 5 for details,
as well as for possible connections to hermitian
$K$-theory which in many ways is the motivic
analog of topological real $K$-theory.

One of the many motivations of this preprint
is that the generalizations of the language
of derived algebraic geometry from
classical to motivic spectra should ultimately
lead to a definition of a motivic
version of $tmf$, generalizing the 
above Theorem \ref{main1} of Lurie about height
2 to the motivic set-up as well.
We will not pursue this in the present 
preprint. (Note that the recent article
\cite{NSO} allows to define motivic elliptic cohomology 
theories and motivic elliptic ring spectra
via motivic Landweber exactness.)
However, concerning motivic derived algebraic geometry,
we wish to point out two interesting applications
of the results of this preprint. First,
motivic symmetric ${\bf P}^1$-spectra
equipped with a suitable positive model
structure satisfy the axioms of 
a $HA$-context of Toen and Vezzosi \cite{TV},
so their machinery applies to this example.
See section \ref{subHA} for a detailed discussion.
Second, the motivic analogue of 
the axioms of Goerss and Hopkins \cite{GH}
is also satisfied (see Theorem \ref{GHaxioms}). I understand that
this second application was established
simultaneously and independently 
by Paul Arne {\O}stv{\ae}r, who wants to use it
for doing motivic obstruction
theory. Both applications are presented 
in section 3.

\medskip

We pause to make some comments concerning the proof of 
Theorem \ref{main2}, which is given in section 4 and uses 
the results established in section 3.
First, one should notice that
the theorem is about $T$-spectra,
but even for $T={\bf P}^1$ the proof involves motivic 
$S^1$-spectra as well. This is mainly due to the fact 
that at some point one needs a motivic version of the recognition
principle which relates $E_{\infty}$-spaces
to connective $S^1$-spectra. The classical
recognition principle is a statement about $S^1$-deloopings, 
and our generalization of it to motivic $S^1$-spectra
is sufficient for our purposes.
Finding a recognition principle for motivic 
${\bf P}^1$-spectra, that is a motivic operad encoding 
${\bf G}_m$- or ${\bf P}^1$-deloopings,
remains one of the main open problems 
in motivic homotopy theory, as already 
pointed out by Voevodsky in \cite[introduction]{Vo2}.
To show that a motivic
version of the recognition principle with respect to $S^1$ 
holds, a previous version of
this preprint invoked
the beautiful ${\bf A}^1$-connectivity theorem of Fabien Morel 
\cite{Mo}, which is known only
over a field. (It seems to be an open question 
if Morel's theorem also holds for 1-dimensional
base schemes. For 2-dimensional
base schemes, there is a counter-example due 
to Ayoub.) However, we later realized that
the proof does not really require this result
and hence holds for more general base
schemes.
\medskip

Keeping in mind these points, the 
strategy for proving
the main theorem may be very roughly
described as follows: First write down
the proof in the classical case, 
that is of Theorem \ref{main1}, choosing
arguments as abstract and as conceptual
as possible. This part is already very interesting 
on its own right, as so far no complete proof of Lurie's 
theorem  in the language of model categories seems to be
available in the literature.
Second, show that this proof 
generalizes to diagram categories (namely
over the Nisnevich site $(Sm/S)_{Nis}$ of smooth varieties
over the given base scheme $S$)
and is well behaved under (left) Bousfield 
localization with respect
to Nisnevich descent and to the 
affine line ${\bf A}^1$. When stabilizing
in the motivic situation, make the right
choice for the circle (sometimes $S^1$ and 
sometimes ${\bf P}^1$ does appear in the proofs)
and for the model structure at every
stage. Indeed, it will turn out that during the
various proofs we have to
consider many different model structures
for symmetric spectra and their lifts
to modules over rings and operads.
Some of these model structures are new, 
and their existence is of independent interest, 
so their presentation here should also serve for
future reference. (Recall that the first 
model structures on motivic symmetric spectra
are due to Jardine \cite{Ja} and Hovey \cite{Ho2}.)
In particular, I am not aware of any discussion
of model structures and derived mapping spaces for 
commutative motivic symmetric ring spectra in the literature so far.

Carrying out the above strategy
requires that various Quillen adjunctions
and equivalences are stable under
localization in a suitable sense,
see e. g. \cite[Theorem 3.3.20]{Hi2} for a result
in this direction. We will provide 
all details in the parts of the proof 
concerning classical spectra. When 
passing to diagram categories and motivic
Bousfield localizations of those, we
will provide details in the first couple
of proofs, but allow ourselves to skip
some of the by then familiar arguments in some of the later proofs.
The reader interested in the classical case
should simply think of the trivial site
and ignore all localization functors
with respect to the Nisnevich topology 
or to the affine line ${\bf A}^1$.
The proof then becomes considerably shorter.
In particular, most (but not all,
see Proposition \ref{Wisgood}) model
structures discussed in section 3 are 
known in that case.

\medskip 

At first glance, it 
might be surprising that we need $E_{\infty}$-structures
to prove a theorem about strictly commutative  
monoids in strictly monoidal model
categories. This is essentially a consequence of
the Lewis paradoxon, as explained at the end of section 2.
If one is only interested in strict adjunctions
and willing to ignore all derived information,
and in particular to sacrify homotopy invariance
of the statement, then there is a much easier proof
not using operads, which we present in section 2.
At the beginning of section 2, we fix some 
notations which will be used throughout this preprint.
 
\medskip

Several results of this article may be generalized
to other left Bousfield localizations of simplicial
presheaves on a site (or even on other diagram categories
satisfying some cardinality conditions), and similarly
for the stable case. See the Remark after Theorem
\ref{anyOmodel} for a more precise statement
in this direction.
 
The referee points out that it should be possible to
prove not only Theorem \ref{main1}, but also
Theorem \ref{main2} using infinity-category techniques.
This would require e.g. to work with the infinity
category of ${\bf A}^1$-local sheaves in some
appropriate infinity-category sense. As I am not
an expert in that field, and as the focus
of this article is on the relevant model structures
for motivic {\it model} rather than {\it infinity}
categories, I will not try to speculate on the details
of this presumably much simpler proof here.

\medskip

This work started as a joint project 
with Niko Naumann, and was presented as such on
a conference in M\"unster in July 2009.
I am indebted to him for the many discussions 
we had on the topics of this preprint. Some parts of
the work presented here have been obtained
in joint work or are at least influenced by these discussions,
and I thank him for allowing me 
to include these parts here. 
Moreover, I wish to thank Stefan Schwede, John Harper and 
Benoit Fresse for discusssions and explanations about certain 
points in their works concerning model structures 
for classical symmetric spectra, operads over them
and $E_{\infty}$-operads, respectively, 
as well as Jacob Lurie for some explanations about \cite{Lu}
and Pablo Pelaez for discussions related to \cite{Pe}
and \cite{Mo}.
 
\section{The non-derived situation: preorientations
which are not homotopy invariant}

The main goal of this section is to establish
the following theorem, which is a non-derived analogue
of Theorems \ref{main1} and \ref{main2}. All the objects involved 
are defined below, and the proof of the theorem is given 
by suitably combining the lemmata in this section.

\begin{theorem}\label{gmadj}
Let $\C$ be an essentially small a category, 
let $M$ be an abelian group object in $\SC$ and 
$A$ an abelian  monoid object in $\SpC$. Then we have a
natural adjunction isomorphism of simplicial sets
$$map_{AbMon \SpC}(S[M],A) \cong map_{AbMon \SC}(M,\Gm'(A)).$$
\end{theorem}

We now introduce some notation. For any monoidal category
$\D=(\D,\otimes)$, we denote the category of monoid objects
in $\D$ by $Mon \D$, and the one of abelian
monoids by $AbMon \D$. When talking about abelian monoids,
we assume moreover that the monoidal categories
and functors involved are symmetric. In our applications, $\otimes$
will be either the cartesian product $\times$
or some smash product $\smash$. All monoids are assumed to be
associative, but not necessarily unital.
We refer the reader to \cite{ML} for precise definitions
of monoidal categories, (strong) monoidal functors etc.

We fix a category $\C$ from now on. In applications $\C$
will be a site, more specifically either the trivial
site or the site $(Sm/k)_{Nis}$ of smooth $k$-schemes, 
$k$ a field, with the Nisnevich topology. 
We denote the category of simplicial presheaves on $\C$
by $\SC$. For a given simplicial presheaf $T$,
we denote the category of presheaves of symmetric $T$-spectra
on $\C$ by $\SpC$. Model structures on these
categories are discussed in \cite{MV}, \cite{Ja}, \cite{Ho2}
and elsewhere, but we won't need them in this subsection.
One might wish to call commutative monoids in $\SpC$ 
``commutative motivic symmetric ring spectra'' in case $\C=(Sm/k)_{Nis}$,
resp. ``commutative symmetric ring spectra'' in case $\C=pt$. \\
Adding a disjoint base point is denoted by $(\ \ )_+$ 
and yields a left adjoint to the functor $F$ forgetting 
the base point in various situations. For a simplicial presheaf
$X$, we set $S[X]:=\Sigma_T^{\infty}(X_+) \in \SpC$
which is defined objectwise as in \cite[Definition 2.2.5]{HSS}
or \cite[Example 1.2.6]{Sc}. 
If $X$ is a monoid in the monoidal category $(\SC,\times)$, then 
$X_+$ is a monoid in the monoidal category $(\SCp,\smash)$
of pointed simplicial presheaves and 
$S[X]$ is a monoid in the monoidal category $(\SpC,\smash)$
of presheaves of symmetric $T$-spectra, see Lemma \ref{adjmon} below. 
We denote the functor sending a presheaf 
of symmetric $T$-spectra to the simplicial
presheaf sitting in degree $0$ by $Ev_0$.

For simplicial presheaves $\F$ and $\G$ we have a simplicial
set $map_{\SC}(\F,\G)$ given by $map_{\SC}(\F,\G)_n:=
Hom_{\SC}(\F \times \Delta^n , \G)$.
We define the simplicial presheaf
$\map_{\SC}(\F,\G)$ by $\map_{\SC}(\F,\G)(c)=map_{\Delta^{op}({\mathcal C}/c)}(\F|c,\G|c)$
where $F|c$ denotes the restriction of the (simplicial)
presheaf $\F$ to the category $\C /c$ of objects in $\C$ lying over
$c$. For presheaves of symmetric $T$-spectra,
we define the simplicial sets $map_{\SpC}$ and 
simplicial presheaves $\map_{\SpC}$ in a similar way. 
Forgetting about simplicial enrichments, 
we write $\Hom$ for the presheaf version of $Hom$.

Finally, we define the (non-derived) multiplicative group
as follows, the derived version of the introduction
being the one using the {\it derived} mapping
space $Rmap$ instead.

\begin{definition}\label{defGm}
The non-derived multiplicative group is
the functor 
$$\Gm':AbMon \SpC \to AbMon \SC$$ 
given by
$$\Gm'(A):=\map_{AbMon \SpC}(S[\Z],A)$$
where $\map_{AbMon \SpC}=
\map_{Mon \SpC}$ is introduced in Definition
\ref{defmapmon} below.
The monoid structure on $\Gm'(A)$ is induced by the 
comonoid structure on $S[\Z]$, the latter lifting
the one on $\Z[\Z]$ corresonding to the
multiplicative group in algebraic geometry
over $Spec(\Z)$.
\end{definition}

\medskip

We define monoidal and strict monoidal functors between
monoidal categories and monoidal transformations
between (strong) monoidal functors as in \cite[chapter XI]{ML}.
An adjunction between monoidal categories is called a {\it monoidal adjunction}
if the unit and the counit are monoidal transformations.
One easily checks that a monoidal functor sends monoids to monoids.
The following Lemma is well-known.

\begin{lemma}\label{adjmon}

(i) We have a monoidal adjunction 
$$(\ )_+ : \SC \stackrel{\rightarrow}{\leftarrow} \SCp : F$$
where $(\ )_+$ is strong monoidal and the forgetful functor
$F$ is monoidal.
Consequently, we have isomorphisms
$$Hom_{Mon \SC}(M,F(N)) \simeq Hom_{Mon \SCp}(M_+,N)$$
for any unpointed monoid $M\in Mon\SC$ and any pointed monoid $N\in Mon\SCp$. 
 
(ii) We have a monoidal adjunction 
$$\Sigma_T^{\infty} : \SCp \stackrel{\rightarrow}{\leftarrow} \SpC : Ev_0$$
where both functors 
are strong monoidal. Consequently, we have isomorphisms
$$Hom_{Mon \SCp}(A,Ev_0(B)) \simeq Hom_{Mon \SpC}(\Sigma_T^{\infty}A,B)$$
for any monoid $A$ in $\SCp$ and any monoid $B$ in 
$\SpC$. 

(iii) We have a monoidal adjunction 
$$S[\ ] : \SC \stackrel{\rightarrow}{\leftarrow} \SpC : F\circ Ev_0$$
where $S[\ ]$ is strong monoidal and $F\circ Ev_0$ is monoidal.
Consequently, we have isomorphisms
$$Hom_{Mon \SC}(M,F\circ Ev_0(B)) \simeq Hom_{Mon \SpC}(S[M],B)$$
for any monoid $M$ in $\SC$ and any monoid $B$ in 
$\SpC$. 

\end{lemma}

\proof
Part (i) is straightforward. The morphisms
of \cite[XI.2.(1),(2)]{ML} for the monoidal
functor $F$ are given by the quotient map $X \times Y \to
X \smash Y$  and by $pt \to (pt)_+=S^0$.
The final statement follows from the obvious remark
that a monoidal adjunction induces an adjunction between
categories of monoids.

Part (ii) is checked objectwise, using \cite[Definition 2.2.5, 
Proposition 2.2.6.1, Definition 2.1.3]{HSS} or
the corresponding results in \cite{Sc} and then proceeding
similar to part (i).

Part (iii) is obtained by 
composing (i) and (ii).
\qed

\begin{definition}\label{defMS}
Let $(\D,\otimes)$ be a monoidal category such that the underlying
category is enriched over 
simplicial sets. We say that $(\D,\otimes)$
satisfies (MS) if there is a natural transformation
of simplicial sets 
$\tau_{x,y,z,w} : map(x,y) \times map(z,w) \to map(x \otimes z,y \otimes w)$
which on $map(\ ,\ )_0=Hom$ coincides with the transformation
sending $(f,g)$ to $f \otimes g$, and we say that 
$(\D,\otimes,\tau)$ is a simplicial monoidal category.
\end{definition}

The property (MS) may be rephrased by saying that
$otimes$ is enriched in simplicial sets.
 
\begin{definition}\label{defmapmon}
Let $(\D,\otimes,\tau)$ be a simplicial monoidal category.
Then for any monoids 
$(x,m_x)$ and $(y,m_y)$
in $\D$, 
we define $map_{Mon}(x,y) \subset map(x,y)$ to be the equalizer of
$map(m_x,y): map(x,y) \to map(x \otimes x,y)$
and $map(x \otimes x,m_y) \circ \tau_{x,y,x,y} \circ
\Delta: map(x,y) \to map(x \otimes x,y)$
where $\tau_{x,y,x,y}$ is as in Definition \ref{defMS}.
If $\D=\SC$, $\D=\SCp$ or $\D=\SpC$, then we denote
the presheaf version of $map_{Mon}$
by $\map_{Mon}$, and the one of $Hom_{Mon}$ by
$\Hom_{Mon}$.
\end{definition}


\begin{lemma}\label{general}
Lemma \ref{adjmon} above remains true when replacing
$Hom$ by $\Hom$ or by $\map$ everywhere.
\end{lemma}
\proof
The isomorphisms for $Hom_{Mon}$ formally imply those
for $\Hom_{Mon}$. The claim about $\map_{Mon}$
follows from Lemma \ref{MS} below.
\qed

\begin{definition}\label{msfunctor}
Let $(\C,\otimes_{\C}, \tau)$ and $(\D, \otimes_{\D}, \tau)$ be simplicial
monoidal categories as in Definition \ref{defMS}.
Let $F:\C \to \D$ be 
a functor of simplicial categories such that
the underlying functor of categories is a monoidal functor
with structure maps $F_2:F(x) \otimes_{\D} F(y) \to F(x \otimes_{\C} y)$
and $F_0:1_{\D} \to F(1_{\C})$. We say that $F$ is a simplicial
monoidal functor if for any objects $x,y$ of $\C$
the diagram
\\
$\xymatrix{
map(Fx,Fy) \times map(Fx, Fy) \ar[d]^{\tau_{\D}} \ar[r]^{F \tau_{\C}} & map(F(x \otimes x),F(y \otimes y)) 
\ar[d]^{F_2^*} \\
map(Fx \otimes Fx,Fy \otimes Fy) \ar[r]_{F_{2*}} & map (Fx \otimes Fx,F(y \otimes y))
}$
\\
commutes.
\end{definition}

\begin{lemma}\label{MS}

(i) Assume that $(\D_i,\otimes_i, \tau_i)$, $i=1,2$ are simplicial
monoidal categories and that 
$$\alpha : \D_1 \stackrel{\leftarrow}{\rightarrow} \D_2 : \beta$$
is a simplicial monoidal adjunction, i. e. $\alpha$
and $\beta$ are simplicial monoidal functors
and there is a monoidal adjunction between
the underlying monoidal functors. 

Then for $(x_i,m_i)$ monoids in $\D_i$, $i=1,2$,
we have an isomorphism 
$$map_{Mon \D_1}(x_1,\beta x_2) \simeq map_{Mon \D_2}(\alpha x_1,x_2)$$
of simplicial sets, and similarly for $\map$.

(ii) The monoidal categories $(\SC , \times), (\SCp , \smash)$ 
and $(\SpC , \smash)$ are enriched over simplicial sets
as categories and satisfy (MS) with respect to the obvious 
choices of $\tau$, hence are simplicial monoidal categories.

(iii) The monoidal adjunctions of Lemma \ref{adjmon} are simplicial.
\end{lemma}
\proof
The proof of (i) is a little long but again straightforward.
In part (ii), for constructing the transformations $\tau$ 
required in (MS) one uses the diagonal $\Delta^n_+ \to \Delta^n_+
\smash \Delta^n_+$, the twist and that for any simplicial
presheaf $K$ (in particular for $K = \Delta^n$) and any $X \in \SpC$ 
one has $K \smash X = (\Sigma_T^{\infty} K) \smash X$
which can be shown objectwise using the results of 
\cite[Chapter I]{Sc}. For part (iii), use that all mapping spaces
involved are defined using the standard cosimplicial object,
and the composition is defined using the diagonal on it. 
\qed

Observe that
for any monoidal category $\C$, one has
$Hom_{Mon \C}(A,B)=Hom_{AbMon \C}(A,B)$ for any abelian monoid objects
$A$ and $B$ in $\C$, and similar for $\Hom$, $map$ and 
$\map$. E. g., both Lemma 
\ref{adjmon} and Lemma \ref{grpmon} below restrict to {\it abelian} monoids.
We now restrict our discussion to unital monoids.
For $M \in Mon \SC$, we denote by $M^{\times}$
the group object in $\SC$ defined by 
$(M^{\times})_k=(M_k)^{\times}$, that is taking 
objectwise the invertible elements in each simplicial degree.
These units satisfy the following.

\begin{lemma}\label{grpmon}

(i) For any $N \in Mon \SC$ a simplicially constant group
object and $M \in Mon \SC$, one has an isomorphism
$$map_{Mon \SC}(N, M) \simeq map_{Groups \SC}(N, M^{\times})$$
and similar for $\map$. In particular, if $N=\Z$
one has $$\map_{Mon \SC}(\Z, M) \simeq M^{\times}.$$

(ii) More generally, if $N$ is a group object in $\SC$, then one has 
an isomorphism 
$$map_{Mon \SC}(N,M) \simeq map_{Groups \SC}(N,M^{\times})$$
and similarly for $\map$.

\end{lemma}

\proof
For simplicially constant $N$, one has isomorphisms \\
$map_{Mon \SC}(N, M)_n \simeq Hom_{Mon \SC}(N, \map_{\SC}(\Delta^n,M))$
$\simeq Hom_{Mon Prshv(\C)}(N,\map_{\SC}(\Delta^n,M)_0) \simeq
Hom_{Mon Prshv(\C)}(N,M_n)$, where the first isomorphism holds
because $( \ \ )_n$ commutes with limits.
The monoid structure on $\map_{\SC}(\Delta^n,M)$ is defined
composing $\tau$, the diagonal $\Delta^n \to \Delta^n \times \Delta^n$
and the monoid structure of $M$.
Using these isomorphisms, part (i) about $map$ reduces to the 
corresponding well-known
result for usual monoids, and the result about
$\map$ follows formally from this by definition.
For the claim about $\Z$ use that
$\map_{Mon \SC}(\Z, M)_n \simeq \Hom_{Mon}(\Z,\map_{\SC}(\Delta^n, M))$
and that $Hom_{Mon}(\Z,M)=M^{\times}$ for usual monoids $M$.
For (ii), one first checks that $\map_{\SC}(\Delta^n \times \Delta^k,
\map_{Mon \SC}(\Z,M)) \simeq \map_{Mon \SC}(\Z, 
\map_{\SC}(\Delta^n \times \Delta^k,M))$ as both are subsets
of $\map_{\SC}(\Delta^n \times \Delta^k \times \Z,M)$
defined by the same diagrams.
As $k$ varies, this implies an isomorphism of simplicial groups
$\map_{\SC}(\Delta^n, M^{\times}) \simeq \map_{\SC}(\Delta^n , M)^{\times}$
where the monoid structure on the right is given by the one on $M$.
Applying $Hom_{Mon \SC}(N,\ )$ and using part (i), one deduces that \\
$Hom_{Mon \SC}(N, \map_{\SC}(\Delta^n, M^{\times})) \\ 
\simeq Hom_{Mon \SC}(N, \map_{\SC}(\Delta^n , M)^{\times}) \\
\simeq Hom_{Mon \SC}(N, \map_{\SC}(\Delta^n , M))$
and the claim now follows by varying $n$, using the adjunction between
$\map$ and $\times$.
\qed
 
Using the above results, Theorem \ref{gmadj} now follows from the 
following chain of isomorphisms using the amplification of the 
indicated results provided by Lemma \ref{general}: 
\[ map_{Mon \SpC}(S[M],A) \stackrel{\ref{adjmon},iii)}{\simeq} map_{Mon \SC}(M,F\circ Ev_0(A)) \]
\[ \stackrel{\ref{grpmon},ii)}{\simeq} map_{Groups \SC}(M,(F\circ Ev_0(A))^{\times}) \] 
\[ \stackrel{\ref{grpmon},i)}{\simeq} map_{Mon \SC}(M,\map_{Mon \SC}(\Z, F\circ Ev_0(A))\]
\[ \stackrel{\ref{adjmon},iii)}{\simeq} map_{Mon \SC}(M,\map_{Mon \SpC}(S[\Z], A)) \]
\[ = map_{Mon \SC}(M,\Gm'(A)).\]
Recall that we may replace $Mon$ by $AbMon$ everywhere.
\qed

\bigskip

Now, what happens if we try to give Theorem \ref{gmadj} 
a homotopy theoretic meaning, equipping everything with
suitable model structures? One problem that may arise
is the definition of
the multiplicative group in Theorems \ref{main1} and \ref{main2}
using the derived mapping spaces with respect to the choosen
model structures. It is not clear if there is a cofibrant replacement of 
$S[\Z]$ which is also a comonoid, that is an affine derived
{\it group} scheme. One may try to show that 
the functor represented by $\Gm'$ is weakly equivalent to
one factoring through simplicial abelian groups.
Independently (in fact maybe not completely independently)
of this, the main problem seems to be the following.
Making our proof homotopy invariant means that all adjunctions 
involved have to be Quillen, and that 
will be impossible to achieve. The problem that appears 
does so already for the trivial category $\C$ with a single object
and no nontrivial automorphisms, that is for classical 
homotopy theory. 
Consider the adjunction of Lemma \ref{adjmon},(iii),
restricted to {\it abelian} monoids.
We want the model structure on $AbMon(\S Sets)$ to be the usual
one. For $\Sp$, we have essentially two families
of model structures, namely the usual ones and the positive 
ones. If we choose a usual non-positive stable model
structure, then this will not lift to a model
structure on $AbMon(\Sp)$ with weak equivalences
and fibrations defined using the forgetful functor
to $\Sp$ because of the Lewis paradoxon, see
e. g. \cite[section 14]{MMSS} or \cite[Remark 4.5]{SS}.
The fact that this adjunction is not Quillen
for any reasonable model structure on abelian monoids 
is why we have to work so much more in the next two
sections, using motivic versions of 
$E_{\infty}$-spectra, of the recognition principle,
of a theorem of Schwede and Shipley \cite{SS2} 
establishing a zig-zag of Quillen equivalences
between $H\Z$-modules in symmetric spectra
and unbounded complexes of abelian groups, etc.
This should not be considered
as a technical problem about model category theory
or symmetric spectra, but as an honest mathematical problem related
to the stable homotopy type of the sphere spectrum and
the content of our main theorem. Therefore it will appear
in some way or another in any language one might choose 
to deal with these questions.

\section{Model structures for algebras over operads in symmetric 
spectra and applications}

The goal of this section is to show that the category of motivic
symmetric spectra as considered by Jardine \cite{Ja}
and Hovey \cite{Ho2} equipped with suitable model
structures satisfies all properties necessary for a motivic
version of derived algebraic geometry.
More precisely, we show (see subsection \ref{subHA})
that motivic symmetric spectra together with 
suitable model structures enable us to construct
model structures for algebras in motivic symmetric
spectra under a given operad, motivic symmetric spectra
satisfiy the assumptions of \cite[section 1.1]{TV}
when choosing suitable additional data
(except that we do not discuss possible
choices of $\C _0$ and $\A$ as introduced in 
\cite[1.1.06 and 1.1.0.11]{TV} here).
They also satisfy a motivic variant of the axioms 
of \cite[1.1 and 1.4]{GH} (see subsection \ref{subGH}).
In particular, we construct model
structures on $\Einf$- and strictly commutative
algebras over motivic symmetric spectra.

At the end of the first subsection, we study model
structures for algebras over motivic operads. This has
not been considered so far. There is also
one new (non-positive!) model structure for $E_{\infty}$-operads
in classical symmetric spectra,
although the existence of such a model structure 
will not be too surprising to the expert.
Namely, Proposition \ref{Wisgood} also applies
to the trivial site $\C$, that is simplicial-set
valued symmetric $S^1$-spectra.
 
Later in this section, some further model structures 
and results related to simplicial presheaves and $H{\bf Z}$-modules
are considered as well.

\medskip

\subsection{Stable model structures}

Let $\C=(Sm/S)_{Nis}$ and fix a cellular left proper model structure on 
$\SC$ which yields the Morel-Voevodsky \cite{MV}
unstable homotopy category $H(S)$, and similarly for the pointed
variant $\SCp$. 
{\it Throughout this section, we will work 
with the motivic injective model structure of} \cite{MV} - or rather with
its extension to simplicial presheaves as in \cite{Ja} - which 
Hirschhorn (\cite{Hi1}, see also \cite[Corollary 1.6]{Hor}) has shown 
to be cellular.
We denote the generating cofibrations (resp. trivial cofibrations)
by $I$ (resp. $J$). Besides being simplicial,
cellular and proper, this model structure has two additional
features which will be important in the sequel.
First, the cofibrations are precisely the monomorphisms,
in particular all objects are cofibrant. Second,
it is a monoidal model category, that is it satisfies
\cite[Definition 4.2.6]{Ho1}. This follows because smashing
with any object preserves weak equivalences, compare
\cite[Lemma 3.2.13]{MV}, \cite[Theorem 1.9]{Hor}
or \cite[Lemma 2.20]{DRO}.
Note that the second condition of loc. cit. for being
a monoidal model category is automatically satisfied
because all objects are cofibrant. 

\medskip

We now fix an object $T$ of $\SCp$. For many arguments below
we may take an arbitrary $T$, but sometimes (e. g. in
Theorem \ref{anyOmodel} and in Theorem \ref{GHaxioms})
we will need that $T \simeq S^1 \smash
T'$ for a suitable $T'$, which holds in particular
for $T=S^1$ and for $T={\bf P}^1$. So we assume that 
$T \simeq S^1 \smash T'$ for a suitable $T'$
from now on, although many results do hold in greater generality.

\medskip

We may apply \cite[Theorem A.9, Definition 8.7]{Ho2} 
and \cite[Theorem 4.1.1]{Hi2} to get a stable model structure on the category
of motivic symmetric spectra $\SpC$ from the above unstable
one on $\SCp$. This model structure
coincides with the one of \cite[Theorem 4.15]{Ja}.
In particular, the {\it motivic stable equivalences} are 
those defined on p. 509 of loc. cit.,
that is defined with respect to injective
stably fibrant objects.

\begin{theorem}(Hovey, Jardine)\label{stableproj}
The above stable model structure on $\SpC$ is simplicial,
proper, cellular and monoidal. 
\end{theorem}
\proof
By \cite[Theorem 4.15]{Ja} we have a proper
closed simplicial model structure. It remains to check the
first condition of \cite[Definition 4.2.6]{Ho1} for a
monoidal model category, that is the pushout-product axiom.
For this we may either apply \cite[Theorem 8.11]{Ho2} as we 
have choosen a model structure on $\SC$ for which all
objects are cofibrant, or directly quote 
\cite[Proposition 4.19]{Ja}.
\qed

This stable model structure on spectra will be referred to as 
the {\it projective stable model structure}. The term
``projective'' refers to the way we obtained the stable
structure from the unstable one, as the unstable
model structure we started with really is an ``injective'' one.
With respect to our fixed choice of the model structure on $\SC$, 
this is a motivic generalization of the model structure considered in 
\cite[Theorem 3.4.4]{HSS}.
It will turn out that this model structure will 
not meet all our requirements, which is
why we need to introduce a motivic version of the (positive) $S$-model=flat
model structure of \cite{HSS} and \cite{Sh}. The reasons 
for considering flat and positive model structures will
become clear below. In the approach of To\"en-Vezzosi,
the reason for considering the flat model structure
is that a motivic generalization of \cite[Corollary 4.3]{Sh}
provides a tool to reduce \cite[Assumption 1.1.0.4 (2)]{TV}
to \cite[Assumption 1.1.0.3]{TV}.

\medskip

We will also need an {\it injective stable model structure}
on motivic symmetric 
spectra, that is a model structure obtained by starting with the levelwise
cofibrations and weak equivalences and then localize to obtain
the stable model structure. 
This is necessary because 
some arguments below will use that the monomorphisms
are cofibrations in a certain model structure, which means that 
for showing that a monomorphism $X \to Y$ is a weak equivalence it is 
sufficient to show that the quotient $Y/X$ is contractible, that 
is weakly equivalent to a point. This model 
structure has been first considered by Jardine \cite{Ja2}.

\begin{theorem}(Jardine)\label{stableinjmodel}
There is a model structure on $\SpC$ with weak equivalences 
being the motivic stable equivalences and
cofibrations being the levelwise monomorphisms. 
This model structure is simplical and
proper.
It is called the
{\it injective stable model structure}.
\end{theorem}
\proof
See \cite[Theorem 10.5]{Ja2} except for right proper,
which follows from the right properness of the stable 
projective model structure which has more
fibrations and the same weak equivalences.
\qed

Next, we establish a flat stable and a positive
flat stable model structure.
The identity morphisms between these four stable 
model structures, that is injective, flat, positive flat 
and projective, are all Quillen equivalences of simplicial
model categories. Of course, it is also possible to 
establish projective and injective positive stable model 
structures, but we will not need these.

As in the classicial case, there is a functor
from symmetric sequences in $\SCp$
to $\SpC$ which is left adjoint to the forgetful
functor. We denote it by $T \otimes -$,
and it enjoys the same formal properties
as the functor $S\otimes -$ in \cite{HSS}.
See e. g. \cite[Definition 2.1.7]{Ho1}
for the definition of $I-cof$ for a set 
of maps $I$ in a category.

\begin{definition}
A map is a {\it motivic flat cofibration}
if it is in $T \otimes M-cof$ where $M$
is the class of levelwise monomorphisms in symmetric sequences.
\end{definition}

\medskip

As we already said above, we will define also define
``positive'' variants of the model structures (at least
for the flat one below),
following \cite[Definition 6.1, Definition 9.1 and p. 484]{MMSS}
which will be necessary to define a model structure on strictly
commutative symmetric ring spectra further below.
This variant has fewer cofibrations than the non-positive
(sometimes also called ``absolute'') model structure.
In particular, the motivic symmetric sphere spectrum 
$\Sigma^{\infty}_TS^0$ is no longer cofibrant,
so the usual contradiction related to the
``Lewis paradoxon'' does not appear
(see e.g. \cite[p. 484]{MMSS}).
Indeed, if one does not work with the positive 
model structure, then in the notations
of Theorem \ref{poscomm} below the condition (2) of 
\cite[Theorem 11.3.2]{Hi2} or equivalently 
\cite[Lemma 2.3.(1)]{SS}
that $U$ takes 
relative $LJ$-cell complexes to stable weak equivalences 
will fail. Looking at the proofs for this condition
(see in particular \cite[Proposition 3.3]{Sh}
and \cite[Lemma 15.5]{MMSS}) one sees how the positive
model structure arises. The key point is that the argument
in the proof of \cite[Lemma 15.5]{MMSS} starting
with ``Since $\Sigma_i$ acts on $O(q)$ as a subgroup
of $O(ni)$'' (read ``on $\Sigma_q$ as a subgroup of
$\Sigma_{ni}$'') does not work if $n=0$. 

\begin{theorem}\label{flatmodel}
The category $\SpC$ abmits a model structure
with the weak equivalences 
being the stable motivic equivalences and cofibrations
being the motivic flat cofibrations. This model
structure is simplicial,
monoidal and proper, and we call it
the {\it flat stable model structure}.
There is also a {\it positive flat stable model structure}
having the same stable weak equivalences
and which enjoys the same properties
(including the pushout product axiom), 
except that the motivic sphere spectrum 
is not flat positive cofibrant.
\end{theorem}
\proof
To establish the flat and the flat positive model structure,
there are various possible proofs. We proceed roughly by
generalizing the corresponding results of \cite{Sh},
see however the remark in parenthesis at the end of the proof.
When adapting Shipley's definition of $I'$ and $J'$ to the motivic case, 
one must work with our above sets $I$ and $J$, of course.
An alternative reference for \cite[Proposition 1.2]{Sh}
which generalizes to the motivic case is \cite{DK}, see
also \cite[Proposition 3.1.9]{Re}. 
The proof of \cite[Proposition 1.3]{Sh} goes through in the
motivic case as well. Note that the model category
on equivariant simplicial presheaves one obtains 
is left proper.
Alternatively, one may deduce the motivic version
of \cite[Proposition 1.3]{Sh} from Proposition 1.2 of loc. cit.
using Hirschorn's \cite{Hi2} or Smith's \cite{Ba} 
generalization of Bousfield localization.

Proposition 2.1 of loc. cit. is just the product
model structure. When establishing the motivic
generalizations of Proposition 2.2 and Lemma 2.3 in \cite{Sh},
the arguments go through and yield a cofibrantly generated level 
model structure, which again is even cellular. 
To see this, observe that (both classical and) motivic
symmetric spectra are cellular because simplicial presheaves
are cellular, and so are products of cellular model
categories. To check the three conditions of
\cite[Definition 12.1.1]{Hi2} for (motivic) symmetric spectra,
first observe that the third condition is
\cite[Proposition A.4]{Ho2}.
We then use the adjunction $(T \otimes - , U=Forget)$
between symmetric sequences and symmetric spectra,
where $U$ commutes with colims.
The proof of the second condition is then
similar to \cite[Lemma A.2]{Ho2}.
Finally, for establishing the first property
one proceeds as in the proof of  
\cite[Proposition A.8]{Ho2}. The argument there
in fact slightly simplifies as we only have to consider
one functor $T \otimes - $ rather than $F_n$
for fixed $n$ with intermediate
considerations concerning $F_m$ for other values
of $m$.
The level flat model structure on motivic symmetric spectra
is left proper because the injective stable model structure
which has more cofibrations and the same weak
equivalences is left proper. To obtain the motivic version of
\cite[Theorem 2.4]{Sh}, that is passing 
from the level to the stable model
structure, one may apply Hirschhorn
localization as in \cite[Definition 8.7]{Ho2}
rather than checking the details corresponding 
to the ones in the proof of \cite[Theorem 2.4]{Sh},
as we have shown that the flat level model
structure is cellular and left proper.
Hence the flat stable model structure is also left proper.
Note that it is right proper because the stable projective model
structure is right proper.
 
The proof that the positive model structure also
exists again goes through in the motivic case.
In more detail, the proof for the positive model structure is 
exactly the same, the only modification being that the motivic
model structure generalizing the one of 
\cite[Proposition 2.1]{Sh} is defined 
as taking on $\Sigma_0$-spaces the cofibrantly
generated model structure
with fibrations and weak equivalences 
being all morphisms. Then take the product
model structure on motivic $\Sigma_n$-spaces
for all $n \geq 0$ as before and proceed as in the 
non-positive case. As the positive model structure has 
fewer cofibrations, it is also left proper.
To show that the positive structure is right proper,
note that the stronger statement
of \cite[Lemma 5.5.3 (2)]{HSS} 
generalizes to the motivic case.
To see this, one uses that
the final argument of loc. cit. carries over
as the ${\bf A}^1$-local model structure
on $\SC$ is right proper by \cite[Theorem A.5]{Ja},
and that the proofs of \cite[Theorem 3.1.14 and Lemma
3.4.15]{HSS} do carry over. 

To check that these model structures are monoidal,
we must check the two conditions of \cite[Definition 4.2.6]{Ho1}.
The second condition in the non-positive case
is easy as the sphere spectrum is 
stably cofibrant because ${*} \to Spec(S)_+$ is cofibrant
in $\SCp$ and $T \otimes -$ is left Quillen.
Hence it remains to check the
first condition, that is the pushout product
axiom. The proof of \cite[Theorem 5.3.7]{HSS}
goes through, and may even be simplified a bit,
see Lemma \ref{flatmonoidal} below, which also applies
to the positive variant.

(It is possible to use the powerful machinery 
of Bousfield localization more systematically
to obtain a different proof. E.g., one may 
quote \cite[Theorem II.4.5]{He}
to obtain the global model structure for $\Sigma_n$-simplicial
presheaves corresponding to \cite[Proposition 1.3]{Sh}
and then impose cardinality bounds to see that
this model structure is cofibrantly
generated and even cellular. This is 
also done in \cite[Theorem 4.9]{Hi1} who attributes
this result to Smith, and in \cite[Theorem 2.16]{Ba}.
Hirschhorn or Smith localization then yields the ${\bf A}^1$-local
model structure on $\Sigma_n$-simplicial presheaves,
with generating sets of (trivial) cofibrations  
different from the ones the approach of \cite{Sh} yields.   
To see that the flat level model structure on $\SpC$
is cellular, one may use the theorem
of Jeff Smith on the existence of left Bousfield
localizations for combinatorial model categories, 
which has been written up recently by Barwick
\cite[Theorem 4.7]{Ba}. Still another variant
would be to apply a more recent localization
theorem of Bousfield as done in the appendix of
\cite{Sc}, and this is certainly not the end of the 
list of variants of proofs...)
\qed

\begin{lemma}\label{flatmonoidal}
Both the flat stable and the positive flat stable model 
structure on $\SpC$ satisfies the push-out product axiom.
\end{lemma}
\proof
We only do the non-positive case. In the positive 
case, the condition on cofibrations follows
similarly, and the one for stable equivalences
then follows from the corresponding property
for the non-positive model structure.

We start by observing that for any finite group $G$,
the above model structure on $G$-objects in $\SCp$
is monoidal because $\SCp$ is monoidal and we have
defined cofibrations and weak equivalences using the 
forgetful functor. It follows that the category of symmetric sequences
in $\SCp$ is monoidal (see e. g. \cite[Theorem 12.2]{Ha2}).
The stable flat model structure on $\SpC$ has 
the same cofibrations as the level structure,
so it remains to show that if given two 
cofibrations $f: A \to B$ and $g:X \to Y$
then if $f$ (or $g$) is a stable equivalence
then so is $f \smash g: A \smash Y \coprod_{A \smash X} B \smash X
\to B \smash Y$.
This can be shown exactly as in 
\cite[Theorem 5.3.7 (5)]{HSS}.
The argument goes through replacing as usual
simplicial sets by simplicial presheaves and
$S^1$ by $T$. In particular
Lemma 3.1.6 of loc. cit. remains valid
in this situation.
\qed

Note that similar to \cite{HSS},
the proof of this lemma provides a variant 
of the above proof that the stable
projective model structure on $\SpC$ 
satisfies the pushout product axiom.
 
\medskip

Next, we wish to study operads $\OO$  
over motivic symmetric spectra.
There are two approaches we are interested in: Simplicial operads,
that is operads in simplicial sets for simplicial monoidal model 
categories, and internal operads in monoidal model categories.
We sometimes apply one and sometimes the other point of view.
Every operad in simplicial sets yields an
internal operad in $\SpC$ via the monoidal functor
$\Sigma_T^{\infty}$, and a similar argument applies to
the unstable case of simplicial presheaves.
The converse is not true, but all operads we are interested 
in are simplical ones. 
We will establish a theorem on the existence 
of model structures for arbitrary internal operads and
stable positive model structures (see Theorem
\ref{anyOmodel}), and weaker 
results in the non-positive case
(see Proposition \ref{Wisgood}).
The latter will be used  when considering
adjunctions of type $(\Sigma^{\infty}_T, Ev_0)$
for $E_{\infty}$-objects.

\medskip

In general, one of the standard ways to construct a model structure
on a category $\D$ is to lift a cofibrantly generated model structure 
on a category $\C$ along a right adjoint in a free/forgetful-style 
adjunction $\C \stackrel{\rightarrow}{\leftarrow}\D$, 
defining fibrations and weak 
equivalences in $\D$ by applying the forgetful functor. 
If this does yield a model structure, then the adjunction is Quillen and 
we say that $\C$ creates a model structure on $\D$.
The main problem when checking the model axioms for $\D$ is 
that in one of the factorizations obtained by the small
object argument it is not clear that certain relative
cell complexes are weak equivalences. See e.g.
\cite[Lemma 2.3]{SS}, \cite[Theorem 11.3.2]{Hi2},
\cite[2.5]{BM1} and certain proofs below.
One strategy for proving this is to establish
a fibrant replacement functor,
see e.g. the discussion in the remark 
after Proposition \ref{levelmodel} below.
Another strategy is to check the required property ``by hand''.
If one is unable to succesfully apply one of
these two strategies, one may - as first done by 
Hovey \cite{Ho0} - try to proceed
by weakening the axioms 
of a model category in a suitable way, which leads to 
the notion of a {\it semi-model category}.
See e.g. \cite{Ho0}, \cite{Sp} and \cite[12.1]{Fr}, we will not pursue this 
approach.

\medskip

The following theorem is a generalization
of a result of Harper \cite{Ha}. Compare also
the article of Elmendorff and Mandell \cite{EM}
which establishes a similar result for
simplicial operads.

\begin{theorem}\label{anyOmodel}
Let $\OO$ be any operad in $\SpC$,
and consider $\SpC$ with the positive flat
stable model structure. Assume
that $T \simeq S^1 \smash
T'$ for some pointed object $T'$.
Then the forgetful functor from
$\OO-alg$ to $\SpC$ equipped with the stable 
flat positive model structure creates a model structure
on $\OO-alg$ in the sense of \cite[Lemma 2.3]{SS}. 
\end{theorem}
\proof
The proof is along the lines 
of \cite[Theorem 1.1]{Ha}. 
We will indicate the nontrivial
modifications to be made in the motivic setting.
One proceeds by showing that the first condition of
\cite[Lemma 2.3]{SS} is satisfied.
Using that transfinite compositions of 
acyclic monomorphisms are acyclic monomorphisms
because they are the acyclic cofibrations for the
model structure of Theorem \ref{stableinjmodel},
this boils down to show the motivic analog of
\cite[Proposition 4.4]{Ha} about a certain morphism $j$.
(Note that in the end we only care about
symmetric sequences of symmetric spectra which are concentrated
in degree zero, as discussed in section 7 of loc. cit..)
To prove the latter, using a filtration argument
it is sufficient (notations taken from \cite{Ha})
to show that $j_t$ is a weak equivalence for all $t$.
For this we proceed as in Proposition 4.29 of loc. cit.,
using that there is a stable model structure on motivic spectra
in which all levelwise monomorphisms are cofibrations,
which exists thanks to Theorem \ref{stableinjmodel}.
Hence we need motivic versions of Proposition 4.28 
and 4.29 of loc. cit.. The proof of Proposition
4.29 uses a five lemma argument  
which requires that smashing with $S^1$ detects stable
weak equivalences, which is fine as $T \simeq T' \wedge S^1$.
Everything else now carries over to the motivic case.
(Note that the positive model structure is used
in \cite[proof of Proposition 4.28]{Ha} 
(``Since $m \geq 1$...''),
and looking at his Calculation 6.15 one sees
exactly what fails for $m=0$.)
\qed

\remark
The proof of Theorem \ref{anyOmodel}
applies to other categories of symmetric
spectra, that is starting with other model structures
on $\SC$ than the motivic one. More precisely,
let us consider the category $\SM$
of simplicial presheaves
on an arbitrary small category $\M$. This category
may be equipped with Heller's global injective 
model structure, which Hirschhorn has shown to be 
cellular and left proper \cite[Theorem 1.4]{Hor}.
Thus we may consider the left (Hirschhorn-)Bousfield
localization $L_S\SM$ with respect to an arbitray set
of morphisms $S$ in $\M$. The proof (see Theorem \ref{Omodelunstable}
below) that this model structure lifts to a model structure on algebras over 
an operad in $\SM$ goes through in this general setting provided we 
can show that the fibrant replacement functor in $\SM$ commutes
with cartesian products. Now to generalize Theorem
\ref{anyOmodel} to symmetric $T$-spectra build from
$\SM$, we need to check that there is an induced
stable injective model structure on symmetric spectra
similar to the one of Jardine (see Theorem \ref{stableinjmodel}),
that is one with all levelwise monomorphisms being 
cofibrations. This together with the fact
that inverting $T \smash - $ inverts $S^1 \smash -$
is all we really need to make the other proofs in this section
work, in particular Theorem \ref{anyOmodel}
and also  Theorem \ref{harperrel} below, applying
the techniques of \cite{HSS}, \cite{Hi2} and\cite{Ha}
in precisely the same way as in the case of motivic
local model structures above.

\medskip

Observe that the above forgetful functor
admits the ``free algebra'' functor as a 
(Quillen) left adjoint, and also that 
both categories admit internal 
mapping spaces compatible under this adjunction.
This will be generalized in Theorem \ref{harperrel} below.

We now shift our attention to non-positive model structures.
The next conjecture is inspired by results
of Harper and Schwede. It applies in particular to 
simplicial $\Einf$-operads. 

\begin{conjecture}\label{Ofreemodel}
Let $\OO$ be any simplicial operad 
such that the action of $\Sigma_n$ on
$\OO(n)$ is free (by which we always mean objectwise
and levelwise free away from the basepoint),
and consider $\SpC$ with the absolute flat or proejective
model structure. Then the forgetful functor from
$\OO-alg(\SpC)$ to $\SpC$ creates a model structure
on $\OO-alg(\SpC)$ in the sense of \cite[Lemma 2.3]{SS}. 
\end{conjecture}
There is the following strategy of proof,
which is an attempt of a motivic generalization
of a variant of a proof for classical symmetric spectra 
as sketched in \cite[section III.4]{Sc}, notations are
again as in \cite{Ha}. In principle, 
it might be applicable to internal operads as well
and not only to those with values in simplicial sets.
As before, we are reduced to consider the 
push-out square of \cite[Proposition 4.4]{Ha}.
It is shown in \cite[Proposition 7.19]{Ha2}
that the $\Sigma_t$-equivariant map
$Q^t_{t-1} \to Y^{\smash t}$ is an acyclic cofibration
of symmetric spectra if $X \to Y$ is an
acyclic cofibration. Using again the motivic generalization
of \cite[Theorem 5.3.7]{HSS}, this implies
that the map  $\OO_A[t] \smash Q^t_{t-1} \to \OO_A[t] 
\smash Y^{\smash t}$ is a monomorphism and a weak equivalence
for any $\OO$-algebra $A$.
Next, we show that the action of $\Sigma_t$
on the motivic symmetric spectra $\OO_A[t]$
is free. By definition, the
action of 
$\Sigma_t$ on $\OO_t(A):=\bigvee_{k \geq 0}\OO(k+t) 
\smash_{\Sigma_k}A ^{\smash k}$ as defined
after \cite[Remark 4.3]{Sc} is free.
Also by definition, $\OO_A[t]$ is an explicit coequalizer of $\OO_t(A)$
with respect to two $\Sigma_t$-equivariant maps coming from another
spectrum with free $\Sigma_t$-action.
Now one has to check that the action of $\Sigma_t$
on this quotient is also free. There are easy examples 
showing that this will not be true for arbitrary operads,
so here one has to use the assumptions on the operad.
Fresse even provided me with an example where
the $\Sigma_n$-action on $\OO(n)$ is free,
but the action on $\OO_n(A)$ is not. 
{\it This is the gap one has to fill when 
using this strategy of proof.}
Hence the diagonal action of $\Sigma_t$ on $\OO_A[t] \wedge Y^{\smash t}$ 
is also free, and so are the ones on 
$\OO_A[t] \smash Q^t_{t-1}$ and on
$\OO_A[t] \smash (Y^{\smash t}/Q^t_{t-1})$.
That is, we have a cofiber 
sequence of $\Sigma_t$-free spaces
$$\OO_A[t] \smash Q^t_{t-1} \to \OO_A[t] \smash Y^{\smash t} \to 
\OO_A[t] \smash (Y^{\smash t}/Q^t_{t-1})$$
for the injective model structure,
which remains true after dividing out the free  
$\Sigma_t$-action on all three objects. 
Indeed, taking coinvariants commutes with taking the cofiber as 
colims commute among each other, and as we have a model structure
in which the monomorphisms are the cofibration,
the cofiber of  $\OO_A[t] \smash Q^t_{t-1} \to \OO_A[t] 
\smash Y^{\smash t}$ is also the homotopy cofiber.
The fact that $\OO_A[t] \smash (Y^{\smash t}/Q^t_{t-1})$ is weakly equivalent 
to a point remains true after dividing out the action
of $\Sigma_t$ as the argument of \cite[Proposition III.4.12]{Sc}
applies to general simplicial monoidal model categories,
including motivic symmetric spectra. Besides the gap above,
this seems to be the only place where it might be an advantage 
to restrict to operads  
defined in simplicial sets rather than in motivic symmetric
T-spectra). Hence 
$\OO_A[t] \smash_{\Sigma_t} Q^t_{t-1} \to \OO_A[t] \smash_{\Sigma_t}
Y^{\smash t}$ is an acyclic cofibration for the 
stable injective model structure, and now the proof can be finished
as the one of Theorem \ref{anyOmodel}.

Note that the recent preprint \cite{GG}
contains a detailed discussion of techniques
related to the problems above.

We can prove 
the above Conjecture \ref{Ofreemodel} at least for the Barratt-Eccles 
operad, which will be sufficient for our purposes.

\begin{definition}\label{BEoperad}
Let $\W$ be the Barratt-Eccles operad with values in the 
monoidal model category $(\S Sets, \times, pt)$,
see e.g. \cite[1.1.5]{BeF} and of course \cite{BE}.
It extends to an internal operad in $\SpC$ via the 
functor $(-)_+$ to pointed simplicial sets,
the constant functor to $\SCp$ and $\Sigma^{\infty}_T$
to $\SpC$ as all of these functors are monoidal. If
we denote the 
internal operad
in $\SpC$ by $\W$ as well, the two
notions of a $\W$-algebra in $\SpC$ thus obtained coincide 
by definition and the above adjunctions.
\end{definition}
  
Observe that in simplicial degree zero, $\W$ is just the 
associative operad $Ass$, and a product of those in higher
simplicial degrees.

We now prove the above conjecture for the operad
$\W$. We don't know of a reference for this result even 
for classical symmetric spectra. The following proof
also applies to other simplicial operads for which
a decomposition pattern similar to the one below for $\W$
may be established.

\begin{proposition}\label{Wisgood}
Conjecture \ref{Ofreemodel} is true
for $\OO = \W$.
\end{proposition}
\proof
According to the above strategy of proof, we must show
that the coequalizer $\W_A[n]$ has a free
action of $\Sigma_n$ for all $n \geq 0$. 
By definition \cite[Construction III.4.8]{Sc},
\cite[Proposition 4.6]{Ha}, the motivic
symmetric spectrum with $\Sigma_n$-action $\W_A[n]$
is the coequalizer of 
$$(m,\W_n(\alpha)):\W_n (\W(A)) \stackrel{\to}{\to} \W_n(A)$$  
where $\W_n(A):= \coprod_{p \geq 0} \W(n+p) \times_{\Sigma_p} A^p$
and $\W(A):=\W_0(A)$.
The map $\alpha: \W(A)=\coprod_{p \geq 0} \W(p) \times_{\Sigma_p} A^p
\to A$ is given by the $\W$-algebra structure of $A$. The map
$m$ is given by the operad structure of $\W$ and 
will be described below, following \cite[Section III.4]{Sc}. 
By definition, the freeness of the $\Sigma_{n}$-action
for a symmetric spectrum has to be checked objectwise (if the site is 
non-trivial, that is in the motivic case), and then for
every symmetric spectrum levelwise, and in each level
degreewise for the simplicial set. Also note that colimits
in motivic symmetric spectra are constructed the same way 
(objectwise, levelwise, simplicially degreewise), and
so are products in simplicial sets and more generally
in simplicial presheaves, and similarly for the
smash product in the pointed case. Finally, the smash
product of a pointed simplicial set with a
symmetric spectrum is given levelwise by the
smash product of the pointed simplicial sets. 
All of this together implies that the whole argument really 
reduces to one of (pointed and even unpointed) simplicial sets, 
so we simplify our notation accordingly.

We now fix $n \geq 0$, and choose orbit decompositions 
of the sets $\Sigma_{n+k}=\W(n+k)_0$ with $\Sigma_n$-action
by left multiplication for all $m \geq 0$. 
These decompositions yield decompositions 
of the simplicial sets $\W(n+k)$ for all $k \geq 0$,
as $W(n+k)$ in simplicial degree $r$ is simply the $r+1$-fold 
product of $\Sigma_{n+k}$ with diagonal
$\Sigma_n$-action, and all simplicial structure maps
are $\Sigma_{n+k}$-equivariant. So we only spell out the 
decompositions in simplicial degree zero.

Recall that as a $\Sigma_n$-set
with action given by left multiplication
the set $\Sigma_{n+1}$ decomposes as a coproduct
of $n+1$ copies of $\Sigma_n$, and inductively
$\Sigma_{n+k}$ decomposes as a coproduct
of $(n+k) \cdot ... \cdot (n+1)$ copies
of $\Sigma_n$ for all $k\geq 0$. We now fix 
particular choices for these decompositions 
of all $\Sigma_{n+k}$ for a fixed given $n$ 
and all $k \geq 0$ once and for all,
and consequently fix decompositions for all
$\W(n+k)$.
 
In every simplicial degree $\W(n+k)$ is a finite set.
The decompositions we choose may be written as products 
$\Sigma_{n+k}=\Sigma_n \times M_{k,n}$
where $\Sigma_n$ acts by left multiplication on the left
factor and trivially on the right factor.
Our decomposition is then determined by the following.
For any positive integer $r$, any element
$\sigma \in \Sigma_r$  is uniquely determined by $\sigma(1,...,r)$.
For $r=n+k$, we write $\sigma = \tau \times \rho$ with
$\tau \in \Sigma_n$ being the element obtained by deleting
all entries in $\sigma(1,...,r)$ which are larger than $n$,
and $\rho \in M_{k,n}$ being determined by where we insert
these remaining elements $n+1, n+2, ..., n+k$ between the 
given permutation of $\tau(1,...,n)$. 

It is obvious that $\W_n(\alpha)$ maps copies of $\Sigma_n$ 
with respect to the above decomposition identically 
(that is not permuting the elements inside each copy
of the $\Sigma_n$-set $\Sigma_n$) to copies 
of $\Sigma_n$. We will show that the same is true 
for the map $m$.
Consequently, the coequalizer $\W_A[n]$ consists of
free $\Sigma_n$-orbits as well, which finishes
the proof.

According to loc. cit., the map 
$m: \coprod_{s \geq 0} \W(n+s) \times_{\Sigma_s} 
(\W(A))^{\times s} \to 
\coprod_{r \geq 0} \W(n+r) \times_{\Sigma_r} A^r$
is defined on the summand for a fixed $s \geq 0$
by the following composition of $\Sigma_n$-equivariant
maps:
$$\W(n+s) \times_{\Sigma_s} \W(A)^{\times s} \stackrel{\cong}{\to} 
\coprod \W(n+s) \times_{\Sigma_s}  (\W(i_1) \times ... \times \W(i_s) 
\times_{\Sigma_{i_1}\times ... \times \Sigma_{i_s}} 
A^{\times i_1 + ...+ i_s}) $$
$$ \stackrel{\cong}{\to} \coprod \W(n+s) \times_{\Sigma_s} 
\W(1) \times ... \times \W(1) \times
\W(i_1) \times ... \times \W(i_s) 
\times_{\Sigma_{i_1}\times ... \times \Sigma_{i_s}} 
A^{\times i_1 + ...+ i_s} $$
$$ \to \coprod \W(n+i_1 + ... + i_s) \times_{\Sigma_{i_1}\times ... \times 
\Sigma_{i_s}} A^{\times i_1 + ...+ i_s}) $$
$$\to \coprod \W(n+r) \times_{\Sigma_r} A^{\times r}$$
where the last map is given by reindexing and the universal
property of coproducts, and the second last map
is given by the structure maps of the Barratt-Eccles operad
$\W$. Now all four morphisms map free $\Sigma_n$-orbits
identically to free $\Sigma_n$-orbits with respect to the 
$\Sigma_n$-decompositions introduced above.
For the first, second and last morphism 
this is obvious. For the third map this can be checked 
using the equivariance condition of the operadic structure 
maps. To see this, again one first looks at the simplicial 
degree zero for which $\W(n+s)_0=\Sigma_{n+s}$, that is the 
associative operad $Ass$. Then generalize to higher degrees 
as explained above, which involves cartesian products of $Ass$, 
and argue componentwise. 
\qed

Once these results are established, one may
deduce the motivic variant of \cite[Theorem 1.4]{Ha}.
Namely, we have the following.

\begin{theorem}\label{harperrel}
Let $f:\OO \to \OO'$ be a morphism of operads
and consider the lifts of the flat 
positive stable model structure on $\SpC$ 
to $\OO-alg(\SpC)$ and $\OO'-alg(\SpC)$.
Then $f$ induces an Quillen adjunction
(enriched over $\S Sets$) 
$$f_* : \OO-alg(\SpC) \stackrel{\rightarrow}{\leftarrow} \OO'-alg(\SpC) : f^*$$
which is a Quillen equivalence if $f$ is a stable 
weak equivalence in every operadic degree.
\end{theorem}
\proof 
By construction of our motivic model structures,
the proof of \cite[Theorem 1.4]{Ha} carries over.
Enrichments are not mentioned in loc. cit., but the arguments given 
there immediately show they behave as well as expected.
\qed

In particular, the model categories of
$\Einf$-algebras and strictly commutative 
monoids in $\SpC$ are Quillen equivalent.

\medskip

We continue to study the absolute (non-positive) situation.
If one can not prove Conjecture \ref{Ofreemodel}
for a given operad $\OO$ using the strategy discussed above, a different
approach might be to first look at the level 
model structure as in the following result. 
\begin{proposition}\label{levelmodel}
There is a projective 
level model structure on $\SpC$ with fibrations and weak equivalences
defined levelwise. This model structure lifts
to $\OO-alg(\SpC)$ for any operad $\OO$. 
\end{proposition}
\proof
In the classical case (i.e. for the trivial site),
the projective level model structure on $\SpC$ is introduced
in \cite{HSS} and the positive variant in
\cite{Sc}, \cite{Sh} and \cite{MMSS}.
These level model structures are cofibrantly
generated with respect to acyclic cofibrations 
$I$ and $J$ (resp. $I^+$ and $J^+$).
To show that it lifts to $\OO-alg(\SpC)$,
as before the only thing that one has to check
is that any map in $\OO J$-cell is a weak equivalence.
For this one may proceed similarly to \cite[Lemma VII.5.6]{EKMM}.
Namely, the geometric realizations of the 
maps in J are ($\Sigma^{\infty}$ of)
inclusions of deformation retracts.
Furthermore, these are stable under
the free functor $\OO$, under push-outs in 
 $\OO-alg(\SpC)$ (by refining an argument
of \cite[Proposition 2.4.9]{Ho1}, as Mandell
kindly explained to me) and under sequential
colimits. Note that geometric realization does preserve colimits,
and a map in  $\SpC$ is a level equivalence
if and only if its geometric realization is. 

In the motivic case (that is for the non-trivial site),
the argument has to be refined a bit. Looking at diagram
categories, one obtains global level (absolute and positive)
model structures on $\SpC$ with $J$ consisting
of inclusions of deformation retracts (objectwise,
in the classical sense) and consequently on $\OO-alg(\SpC)$.
To obtain the ${\bf A}^1$-local level model 
structures, one applies Bousfield-Hirschhorn localization
to $\SpC$ and to $\OO-alg(\SpC)$ to a suitable
set $S$, which yields exactly the cofibrantly 
generated motivic level model structure of \cite[Theorem 8.2]{Ho2} 
on $\SpC$. For $\OO-alg(\SpC)$
one applies Bousfield-Hirschhorn localization with
respect to the set obtained by applying the free
functor $\OO$ to $S$.
\qed

Note that if we could apply Bousfield-Hirschhorn or 
Bousfield-Smith localization with respect to a suitable 
set of stable weak equivalences, this would
yield an alternative proof of Conjecture
\ref{Ofreemodel}. The problem both with Hirschhorn's
and Smith's approach is that $\OO-alg(\SpC)$ with respect
to the level model structure is not left proper
in general, as one sees looking e. g. at $\OO=Comm$.
We will not pursue this approach in our article.

\remark
There is a model structure on the category 
of operads in simplicial sets given by
\cite[Theorem 3.2 and Example 3.3.1]{BM1}, in particular
weak equivalences and fibrations are 
defined on the underlying simplicial sets. 
Note that this model structure 
is different from the one of \cite{Re}.
Looking at  \cite[4.6.4]{BM1} for classical symmetric spectra, 
the argument in the proof
of \cite[Theorem 3.5.(a)]{BM1} applies
to simplicial model categories in general,
defining $\E_f$ in the category of simplicial 
sets and applying $SM 7$.
In particular, we may apply this theorem to
a fibrant replacement map $f:X \to X_{fib}$
in $\SpC$ and a cofibrant model $P$ for the 
$E_{\infty}$-operad with respect to the model structure 
of \cite[Example 3.3.1]{BM1}. Consequently, if we can show 
(which we can't) that the hypothesis in Theorem 3.5 (a) of 
loc. cit. holds, namely
that $f^{\wedge n}$ is a trivial cofibration
for all $n > 0$, then we have constructed a fibrant 
replacement for any object in $\OO-alg(\SpC)$.
Note however that this construction is not functorial.
Hence it is not clear if we may
apply a variant of an argument of Quillen \cite{Qu}, see 
e. g. \cite[Proposition 3.1.5]{Re} or \cite[B.2 and B.3]{Sc0},
to lift the model structure of $\SpC$ to a model
structure on $\OO-alg(\SpC)$. 
In general, in a monoidal model category in which all
objects are cofibrant, the hypothesis holds
as can be shown by an easy induction. More 
generally, to check the hypothesis in a monoidal
model category (a similar argument then presumably
applies to simplicial model categories and simplicial
instead of internal operads), \cite[Remark 3.6]{BM1} claims it 
is enough to have a set of generating trivial cofibrations having 
cofibrant domains. Classical symmetric spectra with the 
(absolute=non-positive) flat or projective model structure satisfy this 
property and are monoidal model categories. But it is not clear why
this helps 
as Harper gave an easy example of a cofibrantly generated
monoidal model category whose generating acyclic cofibrations
have cofibrant domains, but where the above property
for $f^{\wedge n}$ fails already for $n=2$.
In short, the techniques of \cite{BM1} do not provide
an absolute stable model structure for 
$\OO$-algebras in symmetric spectra with
$\OO$ a simplicial or internal cofibrant operad.

\medskip

\subsection{Unstable model structures}

We will now establish a variant
of Theorem \ref{anyOmodel} for spaces
rather than spectra, 
that is an ``unstable model structure'', 
as well as Quillen
adjunctions between unstable and stable
model categories. We also show that a motivic
generalization of the axioms of Goerss-Hopkins 
holds.
For the sake of completeness,
we recall that there is also 
an unstable projective model structure on simplicial
presheaves, starting with fibrations and weak equivalences
defined objectwise and then localizing with respect
to the Nisnevich topology and the affine line
as in the unstable injective case above.
The identity functor obviously induces Quillen equivalences
between the global (and hence between the local)
unstable projective and injective model structure.
Note that the unstable local projective model 
structure is also simplicial, cellular and monoidal
\cite{Bl}, \cite{Hor}. We will only use the unstable injective
model structure in the sequel.

\begin{theorem}\label{Omodelunstable}
For any operad $\OO$,
the injective motivic
model structures on simplicial presheaves lifts
via the forgetful functor to 
the category $\OO-alg(\SCp)$ of $\OO$-algebras
in this model category.
\end{theorem}
\proof
We have a motivic fibrant replacement
functor in $\OO-alg(\SCp)$ as the usual fibrant 
replacement functor commutes with products,
so the argument of Quillen discussed above applies,
see also \cite[Proposition 3.2.5]{Re}.
The existence of the 
motivic (that is ${\bf A}^1$-local) fibrant replacement functor
on  $\OO-alg(\SCp)$ follows from
the existence of a fibrant replacement functor on
$\SCp$ which commutes with $\times$.
More precisely, both motivic fibrant replacement functors
constructed in \cite[Lemma 2.3.20 and Lemma 3.2.6]{MV}
commute with finite limits by \cite[Theorem 2.1.66 and p. 97]{MV}.  
(Note that for the special case
$\OO=\W$, we may of course alternatively adapt the above proof
for $\SpC$.)
\qed

\begin{theorem}\label{easyadj}
There is a Quillen adjunction
$$\Sigma_T^{\infty} : \SCp \stackrel{\rightarrow}{\leftarrow} \SpC : Ev_0$$
where both functors 
are strong monoidal.
Here $\SCp$ is equipped with the above 
injective local model structure,
and $\SpC$ is 
equipped with the projective, flat or injective stable
model structure built from it as discussed above. 
\end{theorem}
\proof
The adjunction was already established in
the first section. By \cite[Proposition 8.5]{Ho2},
the functor $Ev_0$ is right Quillen with respect
to the projective level structure on
$\SpC$, and consequently with respect to all
stable model structures mentioned in the Theorem.
\qed

The following result shows why we it was important
to establish absolute rather than only
positive model structure on $\OO-alg(\SpC)$.

\begin{corollary}\label{adjunstablestable}
Let $\OO$ be an operad such that the forgetful functor to
$\OO-alg(\SCp)$ creates an absolute stable projective resp. flat
model structure on the latter, e.g. $\OO=\W$.
Then the Quillen adjunction of Theorem \ref{easyadj}
induces a Quillen adjuction
$$\Sigma_T^{\infty} : \OO-alg(\SCp) \stackrel{\rightarrow}{\leftarrow} 
\OO-alg(\SpC) : Ev_0$$
for the stable projective resp. flat model structures
on $\OO-alg(\SCp)$.
\end{corollary}
\proof
This follows immediately from Theorem \ref{easyadj}
and the definition of the model structures
on $\OO$-algebras via forgetful functors.
\qed

\subsection{The axioms of Goerss-Hopkins}\label{subGH}

We now show that some of the above model structures
feed into the motivic version of the obstruction
theory machine for $E_{\infty}$ ring
spectra of Goerss and Hopkins \cite{GH}. The machinery 
of To\"en and Vezzosi \cite{TV} will be discussed further below.
As we already said in the introduction, this result
has been obtained independently by {\O}stv{\ae}r.

\begin{theorem}\label{GHaxioms}
Both the flat and the projective stable positive model structures
on $\SpC$ satisfy the motivic analog of the five axioms
in \cite[1.1 and 1.4]{GH}.
\end{theorem}
\proof
Everything has been shown above already except
that the generating cofibrations and the generating
acyclic cofibrations can be choosen to have
cofibrant source and condition (5).
As every object in $\SCp$ is cofibrant, the 
sources of the generating cofibrations 
and of the generating acyclic cofibrations 
for the stable flat model structure 
on $\SpC$ are cofibrant because $T \otimes -$ is
left Quillen, hence this holds in particular for the 
stable flat positive model structure.
Consequently, the same is true for the model structure of operads
as the proof of Theorem \ref{anyOmodel}
shows that the forgetful functor $\OO-alg \to \SpC$
is right Quillen.
To show the claim for the generating acyclic cofibrations, 
note that essentially the same argument goes 
through for those, as the domains of the motivic
variant of the class $K$ of \cite[Definition 4.3.9]{HSS}
also have cofibrant
sources because $Ev_n$ is right Quillen for
the flat level model structure,
hence $F_n$
preserves cofibrations, and the model structures
satisfy the pushout product axiom.  Alternatively,
we may simply quote
\cite[Proposition 4.5.1]{Hi2}.    
Condition (5) follows again from the motivic variants 
of Harper's results discussed above.
Namely, one first uses \cite[Proposition 4.28 (a)]{Ha} 
applied to the cofibration ${*}=Spec(k) \to X$ to see
that $X^{\smash t}$ is then also cofibrant,
and then applies the motivic variant of 
\cite[Proposition 4.29 (b)]{Ha} for 
equivariant symmetric sequences concentrated 
in degree zero, that is equivariant symmetric spectra,
equipped with the model structure of \cite[section 4]{Ha}
for which the weak equivalences are defined
by forgetting the $\Sigma_t$-action. 
\qed

\medskip

\subsection{HA-contexts}\label{subHA}

We now establish model structures for commutative
ring spectra and algebras over those and establish the properties 
required in the axioms of \cite{TV}.
For this, we consider the category $Comm(\SpC)$ of commutative
unital monoids in $\SpC$. 
The notation $Comm$ is taken from \cite{TV},
further below we write $AbMon$ instead which is more consistent with
unstable notations.
The forgetful functor
$U: Comm(\SpC) \to \SpC$ has a left adjoint $L$,
namely the obvious motivic variant of \cite[section 3]{Sh}.
For $R \in Comm(\SpC)$, we define the monoidal category
$R-mod$ in the usual way.

\begin{theorem}\label{modelRmod}\label{Rmodelflatpos}
Consider the stable flat absolute or positive model
structure on $\SpC$, and  
let $R$ be an arbitrary object in $Comm(\SpC)$.
Then there is a model structure
on $R-mod(\SpC)$ where the weak equivalences
and fibrations are defined using the forgetful functor
$R-mod(\SpC) \to \SpC$ where the latter is equipped
with the (absolute or positive) flat model structure.
These model structures are monoidal, proper
and combinatorial. 
Moreover, we have a Quillen adjunctions
$$R \smash - : (\SpC) \stackrel{\rightarrow}{\leftarrow} 
R-mod(\SpC) : U$$
and
$$U : R-mod(\SpC) \stackrel{\rightarrow}{\leftarrow} 
\SpC : Map(R,-)$$
with respect to the flat 
model structure, where we assume that $R$ is cofibrant
for the second Quillen adjunction.
\end{theorem}
\proof
The existence of the model structures follow from the model 
structures of Theorem \ref{flatmodel}
by applying either Kan's lifting theorem \cite[Theorem 11.3.2]{Hi2}
using $R \smash -$ as left adjoint, or the essentially equivalent 
\cite[Theorem 4.1 (2)]{SS}. 
We don't know if the monoid axiom holds
(compare \cite[p. 107]{Ho2}),
but it is sufficient to check
the second condition of \cite[Theorem 11.3.2]{Hi2}
(or equivalently the first condition of \cite[Lemma 2.3]{SS}),
the first one is obvious. If $R$ is cofibrant,
then as the stable flat model structure is monoidal,
the claim follows as explained in \cite[Remark 4.2]{SS}.
In fact, the monoid axiom probably holds in our case by  
some variant of this argument, but we won't need this.

For arbitrary $R$ one must use the motivic generalization
of \cite[Theorem 5.3.7 (5)]{HSS}, compare the classical
proof of \cite[Theorem IV.1.4]{Sc}.
Note that \cite[Assumption 1.1.0.2]{TV} requires the 
model structure also for non-cofibrant $R$.
Both left and right properness in $R-mod$ follow
from left and right properness of $\SpC$. 
Right properness is immediate as the
weak equivalences and fibrations in $R-mod$
are defined by the forgetful functor
the the left proper model category $\SpC$.
To show left properness, one uses that the generating 
cofibrations are level monomorphisms, hence so
are all relative cell complexes built from them.
Now use \cite[Corollary 2.1.15]{Ho2}.

That the two adjunctions exists follows in the
standard way. That the first adjunction 
is Quillen follows from the definition of the model
structures involved.
To prove that the second adjunction is Quillen, one must 
simply show that for cofibrant $R$ 
the functor $Map(R,-)$ preserves fibrations and trivial
fibrations, which immediately
follws from the fact that $\SpC$ is a monoidal
model category.
\qed

It seems possible to prove the above result 
for the projective variant as well, but we won't need this.
Compare also \cite[section 4]{DRO} for similar results
about motivic functors.

\remark
One must show that the above model structure
is combinatorial as desired by To\"en and Vezzosi.
For this, observe that simplicial sets are small \cite[Lemma 3.1.1]{Ho1},
hence so are diagram categories over it (see also
\cite[Remark 2.4]{SS}).
See \cite[Proposition 3.2.13]{HSS} for how to use this to show
that symmetric spectra are also small,
and so are motivic symmetric spectra
using a similar argument.
This shows that the category of symmetric spectra is locally
presentable, and as all model structures we consider
are cofibrantly generated, they are therefore all combinatorial. 
(Note that this also shows that instead
of choosing quite explicit sets of generating (trivial)
cofibrations for the above 
model structures, one might instead take instead
all cofibrations resp. trivial cofibrations with codomains 
bounded dy $\alpha$ (a cardinal that in a suitable sense is
large enough with respect to $\SpC$) as the set 
of generating (trivial) cofibrations.
Then it remains to check that the (trivial) fibrations
are indeed those of the model structure, that is 
it is enough to check the 
lifting property of a (trivial) on these sets.
In most examples this is not hard to see.)
 
\medskip

For $R \in Comm(\SpC)$, we denote the category of
commutative $R$-algebras by $R-Comm(\SpC)$.
We have a forgetful functor
$U:Comm(\SpC) \to \SpC$.

\begin{theorem}\label{poscomm}
The stable flat positive model structure
on $\SpC$ creates a proper combinatorial model structure on
$Comm(\SpC)$ where $f$ is a weak equivalence
(resp. fibration) if and only $Uf$ is.
If $R \in Comm(\SpC)$, then the
same is true for $R-Comm(\SpC)$.
\end{theorem}
\proof 
The existence of the cofibrantly generated
model structure on $Comm(\SpC)$ follows from Theorem 
\ref{anyOmodel} applied to the operad $Comm$.

(In particular, this provides an alternative to the proof
of \cite{Sh} which relies on 
\cite[Theorem 11.3.2]{Hi2},
see also \cite[Theorem A.1.4]{Sc},
\cite[Lemma 2.3]{SS}
and \cite[Proposition 5.13]{MMSS}. 
Note also that \cite{Sc} assumes that $U$ commutes with 
filtered colims which implies that all small colims 
exist, which is an assumption in \cite{Hi2}.
As all objects are small,
the only nontrivial thing of the assumptions that is left 
for Shipley to check is 
that $LJ$-cell complexes - recall that $J$ are the 
generating trivial cofibrations - are stable weak equivalences,
which is not so easy and relies on Propositions
3.3 and 3.4 of her article.)

As discussed above,
it is easy to see that the underlying category is 
locally presentable, and thus the model category is combinatorial.

By the same formal argument as the one
in the proof of \cite[Theorem 3.2]{Sh},
the model structure on $R-Comm(\SpC)$
follows from the one of $Comm(\SpC)$.
It remains to show properness. The model
structures on $Comm(\SpC)$ and more generally 
on $R-Comm(\SpC)$ are right proper 
by the argument of \cite[proof of Proposition 4.7]{Sh}
as we can prove the motivic generalization
of the variant of \cite[Lemma 5.5.3 (2)]{HSS}
for positive level fibrations.
To show that the proof of right properness of 
\cite[Lemma 5.5.3 (2)]{HSS} 
generalizes to the motivic case, one uses that
the final argument carries over
as the ${\bf A}^1$-local model structure
on $\SC$ is right proper by \cite[Theorem A.5]{Ja},
and that the proofs of \cite[Theorem 3.1.14 and Lemma
3.4.15]{HSS} do carry over. 
To check that the model structure is also left
proper, one uses the motivic analogue of 
\cite[Corollary 5.3.10]{HSS} - which 
follows from the motivic generalization of
\cite[Theorem 5.3.7]{HSS} and Ken Brown's lemma - 
and then proceeds as in the proof of
\cite[Proposition 4.7]{Sh}.
\qed

Moreover, we have the following,
which yields the axioms 1.1.0.3 and 1.1.0.4(2)
in the definition of a $HA$-context
as considered by Toen and Vezzosi
\cite{TV}. Observe that axiom 1.1.0.3 is not a formal consequence 
of the property ``monoidal'' established in Theorem \ref{modelRmod}.

\begin{proposition}\label{nolongertodo}
Consider the stable flat positive model structure on
$\SpC$, which we already have shown to be monoidal.

(i) The monoidal model structure on $\SpC$ is symmetric monoidal.

(ii) Let $R \in Comm(\SpC)$. Then for any cofibrant $M \in R-mod$,
the functor $- \smash_R M$ preserves weak equivalences.
For any cofibrant $B \in R-Comm(\SpC)$, the functor
$B \smash _R -: R-mod \to B-mod$ preserves weak equivalences.
\end{proposition}

\proof
Part (i) is clear. 
Part (ii) follows from a motivic generalization of
the ideas of \cite[section 4]{Sh} which can be carried
out thanks to the results we already established.
More precisely, observe that using 
the motivic generalization of
\cite[Corollary 4.3]{Sh} (which holds as 
\cite[Proposition 4.1]{Sh} generalizes to the motivic
situation), Assumption 1.1.0.4 (2) 
reduces to Assumption 1.1.0.3, which is 
a motivic generalization of
\cite[Lemma 5.4.4]{HSS}.
This motivic generalization holds as we have already 
observed that \cite[Theorem 5.3.7, Corollary 5.3.10]{HSS}
generalize to the motivic case. 
\qed

The above results establish
\cite[Assumptions 1.1.0.1 - 1.1.0.4]{TV}
for motivic symmetric spectra $\SpC$.
The non-unital variant (take away the index $0$ in the
definition of the free functor $L$) mentioned in
\cite[1.1.0.4.(1)]{TV} (compare also Remark 1.1.0.5 of loc. cit.)
follows again from Theorem \ref{anyOmodel}
applied to the reduced commutative operad $Comm_{nu}$
with $(Comm_{nu})_0=\Sigma^{\infty}_T(pt)$ when considered
as an internal operad in $\SpC$. 
For classical symmetric spectra, this was already
stated in \cite[Example (4) following Remark 1.1.0.7]{TV}.

\remark
In \cite[p. 20 example (4)]{TV}, it is stated that 
\cite[Assumptions 1.1.0.1 - 1.1.0.4]{TV}
hold for symmetric spectra (that is $\C$ being the trivial
category in our setting) by the results of \cite{Sh}, 
although some of the relevant points are not explained in full detail. Most
of this is carried out in the arguments above.
It remains to show that all categories
involved are locally presentable. The standard references for locally
presentable categories are \cite{AR} and \cite{Bo}.
I do not know a reference for a detailed proof
why symmetric spectra are locally presentable, but 
this easily follws from the smallness property as explained above.
The reason for this condition is that
\cite{TV} quote unpublished work from J. Smith - the relevant
parts are now available thanks to \cite{Ba} - in order to ensure
that certain localized model structures exist, see e. g. 
\cite[section 1.3.1]{TV}. Our localization arguments
in the first half rely on the published work of
\cite{Hi2} on cellular model categories instead,
but the arguments of Smith do apply just as well.
So it is rather a matter of personal taste if one works with
Smith's or with Hirschhorn's version of Bousfield localization.

\medskip

Note that Assumption 1.1.0.4 (2) (which is part of Proposition 
\ref{nolongertodo} (ii)) is probably
not true for the positive projective model structure but only
for the positive flat model structure 
(compare also \cite[Theorem 14.5]{MMSS}). 
Here is why \cite[Proposition 4.1]{Sh} (which is an ingredient of
the proof of \cite[Corollary 4.3]{Sh}) fails 
for the projective model structure already for $R=S$. In the notation 
of loc. cit., the proof uses that the maps of 
${\bf P}_S(S^+I)$ are $S$-cofs, which is deduced from tha fact
that the maps in ${\bf P}(I^{l+})$ are coproducts of monomorphisms
of symmetric sequences. It is not clear that
the corresponding maps ${\bf P}_S(S^+I)$ are stable cofibrations
for the projective model structure.

\medskip

\subsection{$H{\bf Z}$-modules}
Theorem \ref{modelRmod} will be an ingredient in
the proof of Theorems \ref{main1} and \ref{main2}, but
will only be used in the case where 
$T=S^1$ and $R=H\Z$ is the usual
simplicial Eilenberg Mac Lane
spectrum considered as objectwise constant
simplicial presheaf. Such an object $R$ is called
flat resp. projective if it is cofibrant in
$\SpC$ with respect to the flat resp.
projective stable model structure.
So if you don't care about \cite{TV},
the following Lemma allows you to take
a short-cut.

\begin{lemma}\label{HZflat}
Both the classical and the objectwise constant motivic
$S^1$-spectrum $H\Z$
are flat. 
\end{lemma}
\proof
It suffices to show that $H\Z$ is flat
as a classical symmetric spectrum. 
Using the adjunction 
between constant presheaves and presheaves, it
follows that the presheaf of $S^1$-spectra $H\Z$
is cofibrant for the global flat stable
model structure as well,
and  further (motivic left) localizations do not change the 
cofibrations. But as a classical symmetric spectrum, 
$H\Z$ is flat cofibrant, see 
\cite{Sc}.
\qed 

It remains to lift the model structures on ``naive''
$H\Z$-modules and $Ch({\bf Ab})$ to motivic
(=Nisnevich-${\bf A}^1$-local)
model structures on presheaves of those 
as well, using similar techniques as before.
For this we first recall the relevant classical
model structures.

\begin{theorem}\label{abclassical}
The category $Ch({\bf Ab})$ has a model structure
with weak equivalences being the quasi-isomorphisms
and fibrations the epimorphisms.
The full subcategory  $Ch({\bf Ab})_{\geq 0}$
has a model structure with weak equivalences being
the quasi-isomorphisms and fibrations the epimorphisms
in degree $\geq 1$. The inclusion $incl$ and the good 
truncation $\tau_{\geq 0}$ form a Quillen adjunction
between these model categories.
The category $\S \Ab$ has a model structure
with weak equivalences and fibrations 
being the weak equivalences
and fibrations of the underlying simplicial
sets. The Dold-Kan correspondence between
$\S \Ab$ and $Ch(\Ab)_{\geq 0}$ is an 
isomorphism of model categories.
All three model categories are cofibrantly generated and
left proper.
\end{theorem}
\proof
For the cofibrantly generated model
structures, see \cite{Qu} for $Ch({\bf Ab})_{\geq 0}$,
\cite[Theorem 2.3.11]{Ho1} for $Ch(\Ab)$
and \cite[Theorem III.2.8 and Theorem III.2.12]{Ho1}
for $\S \Ab$ (or use the lifting
argument from \cite[Lemma 2.3.(2)]{SS}
as before for the latter, recalling
that there is a fibrant replacement functor
in $\S Sets$ which preserves products).
The claims about the Quillen adjunction resp. 
equivalence are now straightforward.
Left properness for $Ch(\Ab)$, and hence for the other two
model categories as well, follows from 
\cite[Proposition 3.2.9]{Ho1}.
\qed

\begin{theorem}\label{abmotivic}
Theorem \ref{abclassical} generalizes
to the corresponding motivic model categories.
\end{theorem}
\proof
Use the same techniques as before.
First, pass to diagramm categories, that is presheaves
with values in the above model categories, and then 
Bousfield-Hirschhorn-Smith-localize with respect 
to the Nisnevish topology and to the affine line.
Note that $Ch(\Ab)$ is cellular by \cite[Lemma 2.3.2]{Ho1}.
\qed

We will need one more auxiliary model category,
namely ``naive'' $H\Z$-modules,
taken from another article of Schwede and Shipley
\cite[Definition B.1.1]{SS2}. Again, we first explain the classical
case. Once more using the above techniques, everything
generalizes to the motivic case using
\cite[Proposition 12.1.5 and Theorem 13.1.14]{Hi2}
and Bousfield-Hirschhorn localization.
We omit the details.

\begin{definition}
A naive $H\Z$-module is a collection of pointed
simplicial sets $\{M_n\}_{n \geq 0}$
and associative and unital action maps 
$(H\Z)_p \smash M_q \to M_{p+q}$.
A morphism of naive $H\Z$-modules is
a map of graded pointed simplicial sets
which is strictly compatible with the action
of $H\Z$.
\end{definition}

As shown in \cite[Theorem B.1.3]{SS2},
the category $NvH\Z-mod$ of naive $H\Z$-modules
has a model structure 
in which the fibrations and the weak equivalences
are created by the forgetful  
functor $U: NvH\Z-mod \to Sp$ from
naive $H\Z$-modules to classical
Bousfield-Friedlander spectra
with the standard stable model structure 
of \cite[Theorem 2.3]{BF}. The model structure on 
$NvH\Z$-mod is cellular and left proper,
the latter by the same argument as 
in the proof of Theorem \ref{modelRmod}.
One may also consider adjoints to
$U$ as in the case of symmetric spectra
(compare Theorem \ref{modelRmod}), but we won't need
this in the sequel.

\begin{theorem}(Schwede-Shipley)
There is a zig-zag of Quillen equivalences
$$H\Z-mod \stackrel{\leftarrow}{\rightarrow} NvH\Z-mod
\stackrel{\rightarrow}{\leftarrow} Ch(\Ab).$$
\end{theorem}
\proof
See \cite[Appendix B]{SS2}.
\qed

\section{Proof of Theorems \ref{main1} and \ref{main2}}

Having established all necessary model structures
and Quillen adjunctions in the previous section, 
we are now ready to prove Theorems \ref{main1} and \ref{main2}.

\medskip

Let $\M$ be a simplical monoidal
model category and $T$ be a suitable pointed object
in $\M$. We are mostly interested in two cases.
Theorem \ref{main1} is about $\M = \S Sets$
and $T=S^1$.
Theorem \ref{main2} is about $\M$
being simplicial presheaves on $Sm/S$ with the
${\bf A}^1$-local model structure
recalled at the beginning of section 3,
and $T={\bf P}^1$, although throughout
the proof for motivic symmetric spectra
over $T=S^1$ will be considered as well.
Note that the statements of the theorems are
independent of the model structure
one chooses, and by the previous section
there is at least one for all the 
categories involved. Furthermore, we wish to 
emphasize once more that much of the proofs here,
and in fact those of section 3 as well, generalizes 
to other (simplicial) monoidal model categories.
On the other hand, there are some key results which
do not generalize.
In particular, the  Theorems of \cite{SS2}
quoted below, and also Theorem \ref{recognition}
and Theorem \ref{morel}
are really fundamental results specifically
about classical Eilenberg Mac Lance spectra
and delooping along the classical circle
$S^1$, respectively.
Results without $\M$ made explicit
hold in the very general situation
described above. We fix $N \in AbMon(\catM)$ and $A \in AbMon\SpM$,
and assume that $N$ is group-like, 
as defined below.

\bigskip

Putting everything together, we obtain the following diagram 
of categories and adjunctions, with the left adjoint 
displayed on top as usual. For simplicity, we only 
exhibit this diagram in the classical version, that is for 
$\M = \S Sets$. It generalizes to diagram categories,
in particular with the site $\C=(Sm/S)_{Nis}$ as index category,
and to various motivic localizations of those, as  
explained above and below. All categories 
are simplicial model categories, and all adjunctions
are Quillen
(one is even an actual equivalence of categories
preserving the model structure).
The global picture is that there are compatible forgetful 
functors $U$ from the left to the right column.
The functor $V$ in the right column is defined in
\cite[4.3]{HSS} and extends by \cite[B.1]{SS2} to a functor $L$
in the left column. The functor $U$ in the top row
is studied in Theorem \ref{modelRmod}. The upper 
right adjunction is a Quillen equivalence
\cite[Theorem 4.2.5]{HSS} for $\M=\S Sets$ (here
$Sp$ denotes Bousfield-Friedlander spectra \cite{BF}),
and by \cite[Theorem 4.40]{Ja} for 
motivic symmetric $S^1$-spectra.
The left column is explained at 
the end of the previous section.
The dotted adjunction between $E_{\infty}-alg(\S Sets)$
and $Sp$ is standard in other
models for the stable homotopy category, see Lemma \ref{prerecognition}
below. The other dotted arrow to which it restrics
is an equivalence of homotopy categories enriched over 
$Ho(\S Sets)$ as proved in \cite[Theorem 3.45]{ABGHR},
compare our Theorem \ref{recognition} below.
This is a variant of the famous 
{\it recognition principle} due to May \cite{Ma}
and Boardman-Vogt \cite{BV}. We write $E_{\infty}$
for the Barratt-Eccles operad $\W$ introduced in 
Definition \ref{BEoperad}. The precise meaning of the 
dotted arrows in our setting will be explained further below.

\xymatrix{
H\Z-mod \ar[r]^{U} \ar@<-1ex>[d]_U^{\simeq} & \Sp \ar@<-1ex>[d]_U^{\simeq} 
& \\ 
Nv H\Z-mod \ar[r]^{U} \ar@<-1ex>[u]_L \ar@<1ex>[d]^{\H} & 
Sp \ar@<-1ex>[u]_V \ar@{<-->}[dddr] \ar@<-1ex>[d]_{\tau_{\geq 0}}
 & \\
Ch(\Z-mod) \ar@<1ex>[u]_{\simeq} \ar@<-1ex>[d]_{\tau_{\geq 0}} & 
 Sp_{\geq 0} \ar@<-1ex>[u] \ar@{<-->}[dd]^{recogn.pr}_{\simeq}
 & \\
Ch(\Z-mod)_{\geq 0} \ar@<-1ex>[u]_{incl} \ar[d]_{\cong} & & \\
\S \Ab \ar[r]^{U} & E_{\infty}-alg(\S Sets)_{grl} \ar@<+1ex>[r]^{incl} 
& E_{\infty}-alg(\S Sets) \ar@<+1ex>[l]^{GL_1}  
}

\medskip

\begin{lemma}\label{prerecognition}
In the world of Lewis-May-Steinberger
spectra $Sp$ \cite{LMS}, we have a Quillen adjunction of topological
model categories
$$\Sigma^f : E_{\infty}-alg(\S Sets) \stackrel{\rightarrow}{\leftarrow} 
Sp : \Omega^f$$
enriched over topological spaces.
\end{lemma}
\proof
See \cite[Lemma 3.43]{ABGHR}.
\qed

The construction of the functor $\Omega^f$
uses the linear isometries operad which is built in the definition
of Lewis-May-Steinberger spectra, and thus is a key ingredient
when proving the following theorem 
in their setting.

\begin{theorem}\label{recognition}
There is an equivalence of homotopy categories  
$$ Ho(E_{\infty}-alg(\S Sets)_{grl}) \simeq Ho(Sp_{\geq 0}).$$
enriched over $Ho(\S Sets)$.
\end{theorem}
 
We will prove this theorem in our situation,
that is for Bousfield-Friedlander spectra in simplicial
sets $Sp$ and its motivic generalizations.
For this we will use a variant of Lemma \ref{prerecognition} 
which we learned from Schwede, and the symmetric spectra
variant of which will presumably 
be included in the final version of \cite{Sc}.
Namely, we consider the following zig-zag
diagram of Quillen adjunctions

$$Sp \stackrel{\rightarrow}{\leftarrow} Sp (E_{\infty}-alg(\S Sets))
\stackrel{\leftarrow}{\rightarrow} E_{\infty}-alg(\S Sets) $$

\noindent with left adjunctions displayed on top and 
$Sp (E_{\infty}-alg(\S Sets))$ being the category
of spectra with spaces
and spectral structure maps all being $E_{\infty}$.
The left free/forgetful adjunction creates a 
model structure on $Sp (E_{\infty}-alg(\S Sets))$
as usual, that is using the same arguments as
for the existence of the stable model
structure on $Sp$ from $\S Sets$. (This also
can be done in the motivic case below, starting
with the model structure on $(E_{\infty}-alg(\SC)$
established in Theorem \ref{Omodelunstable}.)
In the right Quillen adjunction,
the functor $Ev_0:Sp (E_{\infty}-alg(\S Sets)) \to E_{\infty}-alg(\S Sets)$
is the usual evaluation at the $0th$ space which is 
a right Quillen functor. However its left adjoint is {\it not}
the naive $\Sigma^{\infty}$, but defined using
the simplicial model structure on  $E_{\infty}-alg(\S Sets)$
(which is induced by the Quillen adjunction
with $\S Sets$) when defining the smash products with $S^n$ to define
the level $n$-space of an object in $Sp (E_{\infty}-alg(\S Sets))$.

\medskip

Now by the recognition principle, the left Quillen adjunction
induces an equivalence of (enriched) homotopy categories, and hence
(see e. g. \cite[Proposition 1.3.13]{Ho1}) is an (enriched) 
Quillen equivalence. Moreover, the right Quillen
adjunction induces an equivalence of (enriched)
homotopy categories $ Ho(E_{\infty}-alg(\S Sets)_{grl}) \simeq 
Ho(Sp (E_{\infty}-alg(\S Sets))_{\geq 0})$.
The latter equivalence can also be formulated as
a Quillen equivalence, using suitable localizations 
of the above model structures on $E_{\infty}-alg(\S Sets)$
and $Sp (E_{\infty}-alg(\S Sets)$ as we now explain. 
(We do explain the localized model structure on 
$Sp$ only, for $Sp (E_{\infty}-alg(\S Sets))$, 
the arguments are exactly the same.)
In modern language, see e.g. Hirschhorn \cite[section 5]{Hi2}
or Smith \cite[section 5]{Ba}, these are examples of 
{\it right} Bousfield localizations, that is increasing
the class of weak equivalences while keeping the 
same fibrations. 
I do not know of any published reference
for the following proposition for 
sequential or symmetric spectra, but I learned that 
there is work in progress by Sagave and Schlichtkrull 
- now available, see \cite{SaS} - who
apply similar techniques to study
similar questions for $I$-spaces. 
One should also compare
\cite[section 3.2]{Pe} for a detailed
discussion of how to lift the motivic Postnikov
decomposition to the level of model structures
using right Bousfield localizations, which contains
precisely the arguments needed in our slightly easier case.
Of course, the set $C^0_{eff}$ of loc cit. simply becomes the set
of $S^1$-suspension spectra of the simplicial spheres.

\begin{proposition}\label{rightlocalizedMS}

(i) The category $Sp$ has a simplicial model structure with the 
same fibrations 
as the ones in \cite[Theorem 2.3]{BF}
and weak equivalences the $\pi_n$-isomorphisms
for $n \geq 0$. 

(ii) The category
$E_{\infty}-alg(\S Sets)$ has a simplicial model structure
with the same fibrations as above, that is fibrations
on underlying simplicial sets, and a map
being a weak equivalence if it is one 
after restricting to the invertible
components.
\end{proposition}
\proof
We apply the dual of Bousfield's theorem \cite[Theorems 9.3
and 9.7]{Bou}
to the Postnikov truncation $Q=\tau_{\geq 0}$ 
on $Sp$ resp. to $Q=(-)^{\times}=$ unital components
on $E_{\infty}-alg(\S Sets)$ and the corresponding
transformations $\alpha$. For (ii) the construction
of $Q$ and $\alpha$ is obvious, and for (i) the reader may e.g.
consult \cite[section III.5]{Sc}.
The category $Sp$ is proper by \cite[Theorem 2.3]{BF},
and the properties (A1) and (A2) of
\cite[9.2]{Bou} are obviously satisfied.
The dual of axiom (A3) follows as the hypothesis 
of loc. cit. yield a homotopy pushout square 
and hence isomorphisms $0=Qcone(h) \stackrel{\simeq}
Qcone(k)$ as required. This finishes part (i). For part (ii),
the category $E_{\infty}-alg(\S Sets)$ is right
proper because the forgetful functor to the
proper model category $\S Sets$ preserves limits,
fibrations and weak equivalences.
Showing that it is also left proper is a bit more
subtle, see \cite[Theorem 4]{Sp}
(or \cite[Theorem 12.4.B]{Fr}).
The definition of left proper in loc. cit.
coincides with the usual definition
as all objects in $\S Sets$ are cofibrant.
In order to apply the theorem of loc.
cit. concerning left properness, we need
to now that the operad in question is cofibrant
for the model structure of loc. cit., 
which is created by the one of symmetric sequences. 
That one is equipped with the product model structure
of equivariant simplicial sets which in turn is created
by the one of $\S Sets$ forgetting the group 
action. This implies that the Barratt-Eccles-operad
$\W$ is $\Sigma$-cofibrant, which by definition means
that its underlying symmetric sequence is cofibrant,
as all $\Sigma_n$ act freely. Now we choose a 
cofibrant replacement $\W_{cof}$ of $\W$ (which itself is 
not cofibrant as Fresse kindly explained to me),
which then in particular is also
$\Sigma$-cofibrant (see e. g. \cite[Proposition 4.3]{BM1}).
Finally, the model categories
$\W-alg$ and $\W_{cof}-alg$ are Quillen equivalent
by \cite[Theorem 12.5.A]{Fr},
so we do not distinguish between them in our notations
in the sequel.
The dual statements of the conditions (A1), (A2) and (A3)
of \cite{Bou} are again easy to check.
\qed

Putting everything together and again suppressing the Quillen
equivalence between 
$Sp$ and $Sp (E_{\infty}-alg(\S Sets))$ in our notations,
we obtain the following square of Quillen adjunctions
(model structures omitted from the notations). It 
corresponds to the lower right triangle in 
the large diagram of Quillen adjunctions above before
Lemma \ref{prerecognition}.

\xymatrix{
E_{\infty}-alg(\S Sets)
\ar@<+1ex>[r]^{Id} \ar@<+1ex>[d]^{H}_{\simeq} & E_{\infty}-alg(\S Sets)
\ar@<+1ex>[d]^{H} \ar@<+1ex>[l]^{Id} \\ 
Sp \ar@<+1ex>[r]^{Id} \ar@<+1ex>[u]^{Ev_0} & 
Sp \ar@<+1ex>[u]^{Ev_0} \ar@<+1ex>[l]^{Id} 
}

We claim that the left vertical pair
is a Quillen equivalence. The horizontal equivalences
are the right Bousfield localizations we just described.
In more detail, the cofibrant objects
in the left hand side model categories are the cofibrant objects
with respect to the model structures on the right hand side 
which are moreover group-like $E_{\infty}$-spaces resp. $(-1)$-connected
spectra. Hence the Quillen adjunction on the
right hand side induces one on the left hand side.
This is an example of \cite[Theorem 3.3.20(2)(a)]{Hi2}:
note that \cite[Definition 8.5.11 (2)(a)]{Hi2} applied
to right localization with respect to $\tau_{\geq 0}$ does not
produce additional weak equivalences in $E_{\infty}$-algebras
and therefore induces a Quillen adjunction between the lower left and 
the upper right corner in the above diagram. Composing this Quillen adjunction
with the one on top leads to the one of the left hand side
we are looking for.
As we already pointed out above, this Quillen adjunction then induces an
equivalence of homotopy categories by the recognition priciple, 
hence it is a Quillen equivalence (see e. g. \cite[Proposition 1.3.13]{Ho1}).

Recall that when writing
$E_{\infty}-alg(\S Sets)_{grl}$ resp. $Sp_{\geq 0}$
in the large diagram above
and further below, we really mean the model
categories $E_{\infty}-alg(\S Sets)$ resp. $Sp$
with the right localized model structures
established in Proposition \ref{rightlocalizedMS}.
The total right derived functors of the right
adjoint identity functors in the above square are precisely
$GL_1$ resp the Postnikov functor $\tau_{\geq 0}$,
thus justifying the labels on the arrows in the big diagram
further above. 
This finishes our discussion
of the proof of Theorem \ref{recognition}
and its refined formulation in the language of model  
categories.
\qed

\medskip

Again, all above Quillen adjunctions and equivalences
in the above diagram generalize to the motivic situation 
using always the same kind of arguments involving 
left Bousfield localizations of diagram categories.

\begin{theorem}\label{morel}
The above Quillen adjunctions 
induce Quillen adjunctions for the corresponding motivic
categories, which then induce an equivalence
of homotopy categories 
$$ Ho(E_{\infty}-alg(\SC)_{grl}) \simeq Ho(Sp^{S^1}(\C)_{\geq 0})$$
enriched over $Ho(\S Sets)$ and even over
$Ho(\SC)$.
\end{theorem}
\proof
The above Quillen adjunction between $E_\infty$-spaces and spectra
generalizes to one between global model structures
on diagram categories (see e. g. \cite[Theorems 11.6.1
and 11.6.5]{Hi2}). We claim that this one then induces a Quillen
adjunction after left Bousfield localization on both sides
with respect to the Nisnevich topology and to ${\bf A}^1={\bf A}^1_S \to S$ 
by standard arguments, that is applying \cite[Theorem 3.3.20 (1)(a)]{Hi2}.
On $(E_{\infty}-alg(\SC))$, this is precisely 
the model structure established in Theorem \ref{Omodelunstable}.
Indeed, looking at the fibrant replacement functors discussed there,
we see that we obtain the correct ${\bf L}F\C$ in the notation of loc. cit..
As before, we then wish to apply a suitable right
Bousfield localization to 
these left localized motivic model structures, thus  
obtaining the homotopy categories 
of connected motivic $S^1$-spectra 
and grouplike motivic $E_{\infty}$-spaces 
using the (diagram versions of)
the right Bousfield localizations considered in Proposition
\ref{rightlocalizedMS}.
To see that the right localization of the ${\bf A}^1$-local
structures on motivic 
$E_{\infty}$-spaces exists, we may apply the dual of 
\cite{Bou} as before.
The arguments above imply that motivic $E_{\infty}$-spaces are 
cellular and left proper. To see that they are right proper,
recall that motivic spaces are right proper
and the ${\bf A}^1$-local model structure on
$E_{\infty}$-spaces is created by the forgetful functor
which preserves pull-backs.
Then one checks the dual of the remaining hypotheses
of \cite{Bou} with respect to $Q=(-)^{\times}$.
Concerning the right localization of motivic $S^1$-spectra, it turns out 
to be more convenient to apply
\cite[Theorem 5.1.1]{Hi2} rather than \cite{Bou}.
That is, we proceed as Pelaez does in \cite[section 3.2]{Pe}.
Recall that motivic $S^1$-spectra are cellular and proper
by Hovey and Jardine, that is (the sequential
spectra version of) Theorem \ref{stableproj}.
From the (sequential spectra version of) 
\cite[Proposition 3.2.4]{Pe} it easily follows
that the right Quillen functor $Ev_0$ 
with respect to the motivic model structures
remains right Quillen when applied to the right localizations
we just described.
As it induces an equivalence of homotopy
categories, it is a Quillen equivalence as
claimed.
\qed

As Pelaez explained to me, Morel's connectivity result
(which is valid only for $S=Spec(k)$) really is stronger
than what we have used here. It can of course not be recovered
using only the above techniques.

\bigskip

Now we are ready for the proof of the main
theorems. We write down a chain 
of natural weak equivalences of simplicial sets,
so applying $\pi_0$ yields Theorems \ref{main1} and \ref{main2}.
To simplify notation, we drop all base points in the sequel.
 
Let $\catM$ be one of the two monoidal model categories
we are interested in, see the beginning of this section.
We have

$Rmap_{AbMon(\SpM))}(\Sigma^{\infty}_{T}N,A)$ 

$\simeq Rmap_{\Einf(\SpM))}(\Sigma^{\infty}_{T}N,A)$

using Theorem \ref{harperrel} and flat positive
stable model structures.  Now 
the identity is a Quillen equivalence
between the positive and the non-positive 
(see Theorem \ref{anyOmodel} and Proposition \ref{Wisgood}) 
model structure, so we may switch
to the latter on $\Einf(\SpM))$.
Then using Corollary \ref{adjunstablestable}, 
we have

$\simeq Rmap_{\Einf(\catM)}(N,Ev_0(A))$

As $N$ is an abelian
group by assumption, hence grouplike, we have

$\simeq Rmap_{\Einf(\catM)_{grl}}(N,GL_1(Ev_0(A)))$ 

by Proposition \ref{rightlocalizedMS} above. Recall the meaning 
of the heuristic notations $\Einf(\catM)_{grl}$ and $GL_1$ as introduced 
immediately after loc. cit., it would be more accurate to say
that ``$N$ is cofibrant in the right localized model
structure of Proposition \ref{rightlocalizedMS}''. 
The chain of weak equivalences continues with

$\simeq Rmap_{Sp_{\geq0}}(HN,gl_1(A))$

using Theorem \ref{recognition} resp. Theorem \ref{morel}
and its proof, which 
defines $HN$ and $gl_1(A)$. Observe that 
$GL_1 \circ Ev_0=Ev_0 \circ gl_1$, and moreover $gl_1$
allows a model-theoretic description as right Quillen adjoint
similar to $GL_1$. As before, we have suppressed the left hand 
side Quillen equivalence in the zig-zag of the Quillen adjunctions
after Theorem \ref{recognition} from our notations.
Recall also that we are dealing with $S^1$-spectra and the 
usual Eilenberg-Mac Lance spaces here.
The next weak equivalence

$\simeq Rmap_{Nv H\Z-Mod(Sp_{\geq0}(\catM))}(HN,Rmap_{Sp_{\geq0}(\catM)}
(H\Z,gl_1(A)))$

is just a formal adjunction, see Theorem \ref{modelRmod}
which allows a variant for
``naive'' $H\Z$-modules.
Note in particular that $HN$ is a module over $H\Z$.


$\simeq Rmap_{Ab(\catM)}(N,Rmap_{\Einf(\catM)_{grl}}(\Z,GL_1(Ev_0(A))))$

using Theorem 3.23 of Schwede-Shipley resp. its motivic generalization
and Theorems \ref{abclassical} and \ref{abmotivic}. 
Note that this is compatible with Theorem \ref{recognition}
by the large commutative diagram above.
In particular, $Rmap_{\Einf(\catM)_{grl}}(\Z,GL_1(Ev_0(A)))$
no longer denotes a (presheaf of)
$H\Z$-module(s), but the corresponding
(presheaf of) simplicial abelian group(s).
Observe that the model structures on $Ab(\catM)$ and
$AbMon(\catM)$ are compatible since both are created via
the forgetful functor to $\catM$.

$\simeq Rmap_{AbMon(\catM)}(N,Rmap_{\Einf(\catM)}(\Z,GL_1(Ev_0(A))))$

$\simeq Rmap_{AbMon(\catM)}(N,Rmap_{\Einf(\catM)}(\Z,Ev_0(A)))$

as $\Z$ is grouplike again by Proposition \ref{rightlocalizedMS} above

$\simeq Rmap_{AbMon(\catM)}(N,Rmap_{\Einf(\SpM)}(\Sigma^{\infty}_{T}\Z,
A))$

and finally proceeding as above

$\simeq Rmap_{AbMon(\catM)}(N,Rmap_{AbMon(\SpM)}(\Sigma^{\infty}_{T}\Z,
A))$.

\medskip

Note that the above chain of weak
equivalences really arises from ``zig-zags'', as various of the 
enriched Quillen equivalences
in the above argument go in the ``wrong'' direction.
To start with, the identity is a left Quillen adjoint from 
the positive to the absolute model structure
on symmetric spectra and algebras over those.
When considering derived mapping spaces, 
we must choose a cofibrant and a fibrant replacement 
functor for both model structures. In this situation,
we may simply choose 
the cofibrant replacement functor with respect to
positive model structure (which then also is one for
the absolute model structure) and the fibrant
replacement functor with respect to the absolute model
structure (which then is also one for the positive model
structure). These choices show that the chain of weak
equivalences leading to the Main Theorems really can be 
choosen to be one with is natural both in $N$ and $A$.
Another such zig-zag is hidden in the recognition principle,
and still another one in the motivic generalization of \cite{SS2}.

\medskip

Putting everything together, we therefore have
a natural weak equivalence of derived simplicial
mapping spaces

$$Rmap_{AbMon(\SpM))}(\Sigma^{\infty}_{T}N,A)$$ 
$$\simeq Rmap_{AbMon(\catM)}(N,Rmap_{AbMon(\SpM)}(\Sigma^{\infty}_{T}\Z,
A))$$ which after aplying $\pi_0$ finishes the proof of 
Theorems \ref{main1} and \ref{main2}.

\section{Motivic preorientations and orientations of
the derived multiplicative group}

In this section, we explain how the Theorems
\ref{main1} and \ref{main2} lead to the results 
about orientations and $K$-theory stated in the introduction.

\medskip

Theorem \ref{main2} applies to $N=\overline{W}(GL_1)\simeq{\bf P}^{\infty}\in\SC$. Here, $\overline{W}(\cdot)$ is a specific model for the classifying spave of a simplicial group, see \cite[Chapter V, 4]{GJ}. It is easy to see that $\overline{W}(\cdot)$ sends commutative simplicial abelian groups to commutative monoids in $\SC$. The equivalence $\overline{W}(GL_1)\simeq{\bf P}^{\infty}$
is a special case of  \cite[Proposition 3.7]{MV}.
Beware of the difference between $GL_1$ and the $\Gm$ above.
In Theorem \ref{main1}, the same argument applies to
the topological group $S^1=U(1)$
with classifying space ${\bf CP}^{\infty}$. 
 
We do not suggest a definition of the notion of a derived group 
scheme in $AbMon \SpC$ here.
Compare \cite[section 3]{Lu} for a motivation of
the following definition, at least for the trivial site.

\begin{definition}

A pre-orientation on a derived group scheme $G$ 
over a motivic symmetric spectrum $A$ is an element in $Hom_{AbMon \SC}(\overline{W}(GL_1),G(A))$.

\end{definition}

Note that this definition is related, but not
equivalent to more classical notions of orientations
as e. g. in Adams book \cite{Ad}.
(Namely, in contrast to \cite{Ad}, we consider strict monoid
homomorphisms to the infinite dimensional
projective space, and also we do consider such maps only up to
homotopy.)
However, by \cite{Lu} it is the ``correct'' definition
in order to obtain the right definition of $tmf$,
and in the height 1 case to obtain $KO$.

\medskip

One easily checks that for any simplicial abelian monoid,
the associated suspension spectrum
is a commutative motivic ring spectrum
(compare \cite[Example I.2.32]{Sc} and Lemma \ref{adjmon}).
Theorem \ref{main2} may thus be rephrased as follows
in the special case of ${\bf P}^{\infty}$, where
we write $S[-]$ for $\Sigma^{\infty}_T(-)$
following Lurie's notation (as introduced in the beginning 
of section 2 already).

\begin{theorem}\label{premain}
There is a bijection between preorientations of
$\Gm$ over $A$ and $Hom_{Ho(AbMon \SpC)}(S[\overline{W}(GL_1)],A)$.
In other words, $S[\overline{W}(GL_1)]$ classifies preorientations 
of the derived multiplicative group.
\end{theorem}

Lurie also gives a definition of an orientation,
see \cite{Lu}. We do not suggest a motivic generalization
of this definition in general, either.
However, any reasonable generalization of the notion
of an orientation from the classical to the
motivic case will certainly imply a bijection between the set 
of orientations on $\Gm$ over $A$ and the set 
$Hom_{Ho(AbMon \SpC)}(S[\overline{W}(GL_1)][\beta^{-1}],A)$
for a certain lift of the motivic Bott element $\beta$ (see below).
In other words, $S[\overline{W}(GL_1)][\beta^{-1}]$ 
will classify orientations 
of the derived multiplicative group.

\medskip

By recent work of Spitzweck-{\O}stv{\ae}r \cite{SO}
and independently of Gepner-Snaith \cite{GS},
we have the following algebraic version
of Snaith's theorem.

\begin{theorem}\label{algsnaith}
(Spitzweck-{\O}stv{\ae}r, Gepner-Snaith)
There is an isomorpism of commutative monoids in $SH(S)$ between 
the underlying motivic
spectrum of $S[\overline{W}(GL_1)][\beta^{-1}]$
and Voevodsky's motivic spectrum representing algebraic
$K$-theory \cite[section 6.2]{Vo} where $\beta \in
\pi_{2,1}(BGL_1)$ is a lift of the motivic Bott element.
\end{theorem}

In light of this result, 
our Theorem \ref{premain} above
may be rephrased by saying that 
{\it algebraic $K$-theory classifies 
orientations of the derived multiplicative group}.
More precisely, one must either assume $S$ regular
here or work with Weibel's homotopy invariant
algebraic $K$-theory \cite{We} for non-regular base schemes.

The classical Snaith theorem \cite{Sn} together with
some considerations about suspension spectra
and suitable localizations of those being semistable
symmetric ring spectra leads to a description
of topological $K$-theory as a strictly commutative ring
spectrum, that is an abelian monoid in symmetric spectra.
More precisely, suspension spectra
are ``semistable'' in the sense of \cite[Theorem I.4.42]{Sc}
by \cite[Example I.4.46]{Sc}. It follows that
$S[\overline{W}(GL_1)][\beta^{-1}]$
is again a symmetric spectrum by 
\cite[Corollary I.4.67]{Sc}, hence
complex topological $K$-theory is represented
by a strictly commutative ring spectrum in 
$\Sp$. This argument generalizes to the motivic situation.

\begin{proposition}(R\"ondigs,Spitzweck,{\O}stv{\ae}r)
The object $S[\overline{W}(GL_1)][\beta^{-1}]$
is a commutative monoid in $\SpC$.
\end{proposition}
\proof
See \cite{RSO}.
\qed

Returning again to the classical case, 
the fact that $K_{top}^{h\Z/2} \simeq KO_{top}$
implies that $KO$ classifies all oriented
derived multiplicative groups, see \cite[Remark 3.12]{Lu}
for details. A similar statement for algebraic and hermitian $K$-theory,
namely $K_{alg}^{h\Z/2} \simeq KO_{alg}$
was conjectured to hold for arbitrary rings
with 2 invertible at least after a suitable completion,
see \cite[3.4.2]{Wi}. This conjecture has been proved
in many cases, see \cite{Ko}, \cite{BKO} and more recently
\cite{HKO}, \cite{BKSO}, but in 
general it is wrong as \cite{BKO} show. 

\medskip

Finally, let us mention that there are of course many examples
of abelian monoids in symmetric $T$-spectra. 
Suspension spectra of abelian monoids, e.g. of algebraic
groups or abelian varieties, are obvious examples.
Another example is Voevodsky's  algebraic
cobordism spectrum ${\bf MGL}$, as explained in 
\cite[section 6.5]{PY}, \cite[section 2.1]{PPR}. The techniques of 
Schlichtkrull \cite{Sk} then yield many 
more examples, as the proof of Theorem 1.1
of loc. cit. carries over to the motivic Thom spectrum,
and hence applies to $\I \U / BGL$ with $\U$
being the category of motivic spaces, that is
simplicial presheaves. Moreover, one may 
try to use the isomorphism ${\bf A}^n-{0} \simeq 
S^{n-1} \smash \Gm ^n$
to extend the picture to a motivic version of generalized 
Thom spectra with respect to a motivic version
of $BF$, that is with self-maps on $T^n$.
We might pursue this topic in some other article.

\newpage

\section*{Appendix: Erratum from 2017}

\setcounter{theorem}{0}
\renewcommand{\thesection}{A}
\def\A{\mathbf{A}}

\begin{abstract}
We correct a claim concerning
motivic $S^1$-deloopings. 
\end{abstract}

\medskip

In \cite[Theorem 4.4]{Ho}, we claim that 
the classical recognition principle of May
et al carries over to motivic $S^1$-spectra
with respect to the $\A^1$-local model structure.
Unfortunately, the proof of this theorem is incomplete.

Moreover, the article \cite{Ch} shows that
a certain $\A^1$-local sheaf $\Z(\Gm)$
of abelian groups
is not strongly $\A^1$-invariant. Hence
this sheaf, considered as a group-like object 
in $\SC$, has only $S^1$-deloopings which
are not $A^1$-local, contradicting Theorem 4.4
of loc. cit. The correct version of Theorem 4.4 is stated
and proved in the recent preprint \cite{EHKSY}
of Elmanto, Hoyois, Khan, Sosnilo and Yakerson.
In short, one has to take care of the
interplay between group completion and 
$\A^1$-localization. If we just consider the
$Nis$-local model structure before $\A^1$-localization,
Theorem 4.4 is true by \cite[Theorem 5.2.6.15]{Lu}. 
  
\medskip

Consequently, the proof of Theorem 1.2 of
\cite{Ho} is incomplete for the $\A^1$-local structure
and complete only $Nis$-locally. Currently we do not know
if Theorem 1.2 is true as stated; hopefully future research will
answer this question. On the other hand,
the mistake does not affect any of the results of section
3 of \cite{Ho}, where various model structures related
to motivic operads are studied.

\medskip

Acknowledgement: We thank Niko Naumann
for asking questions about \cite{Ho}
and pointing out the references \cite{Ch} and \cite{Lu},
Utsav Choudhury for answering several
questions related to \cite{Ch} and the authors of \cite{EHKSY}
for sharing a preliminary version of their work 
with me.

Jens Hornbostel, 
Bergische Universit\"at Wuppertal,
FB C, Mathematik und Informatik,
Gau{\ss}strasse 20, 42119 Wuppertal,
hornbostel@math.uni-wuppertal.de.

\end{document}